\documentclass[12pt]{article}
\usepackage{lmodern}
\usepackage{lmodern}
\usepackage[T1]{fontenc}
\usepackage{textcomp}
\usepackage[utf8]{inputenc}
\usepackage{xcolor}
\usepackage{verbatim}
\usepackage{mathrsfs}
\usepackage{mathtools}
\usepackage{enumitem}
\usepackage{amsmath}
\usepackage{amsthm}
\usepackage{amssymb}
\usepackage[a4paper]{geometry}
\usepackage{microtype}
\usepackage[bookmarks=false,
 breaklinks=false,pdfborder={0 0 1},backref=false,colorlinks=true]
 {hyperref}
\hypersetup{
 linkcolor=blue,citecolor=blue,urlcolor=blue}

\makeatletter

\usepackage[utf8]{inputenc}
\usepackage{amsmath}
\usepackage{amssymb}
\usepackage{amsthm}
\usepackage{mathtools}
\usepackage{graphicx}
\usepackage{tikz-cd}
\usepackage{mlmodern}
\usepackage{fouriernc}
\usepackage[T1]{fontenc}
\usepackage{xcolor}
\usepackage{pdfrender}
\usepackage{mathtools}        
\usepackage{bm}               
\usepackage{mathrsfs}         
\usepackage{dsfont}           
\usepackage{cancel}           
\usepackage{graphicx}
\usepackage{booktabs}

\usepackage{amsthm}\usepackage{mathrsfs}
\usepackage{bookmark}
\numberwithin{equation}{section}
\newtheorem{theorem}{Theorem}[section]
\newtheorem{lemma}[theorem]{Lemma}\newtheorem{proposition}[theorem]{Proposition}\newtheorem{corollary}[theorem]{Corollary}\theoremstyle{definition}
\newtheorem{example}[theorem]{Example}\theoremstyle{remark}
\newtheorem{remark}[theorem]{Remark}\newcommand{\D}{\mathbb D}
\newcommand{\T}{\mathbb T}
\newcommand{\R}{\mathbb R}
\newcommand{\C}{\mathbb C}
\newcommand{\E}{\mathbb E}
\DeclareMathOperator{\Var}{Var}
\DeclareMathOperator{\Law}{Law}

\DeclareMathOperator{\Area}{Area}
\title{Symmetrization and Optimal Planar\ Skorokhod Embedding}
\author{Maher Boudabra\thanks{\raggedright Department of Mathematics, Monastir Preparatory Engineering Institute, Tunisia. \mbox{Email: \href{mailto:maher.boudabra@gmail.com}{maher.boudabra@gmail.com}.}}}
\date{Revised working manuscript with ratio theorems\\12 September 2026}

\makeatother

\providecommand\theoremname{Theorem}
\theoremstyle{plain}

\begin{document}
\title{Quantitative Brownian Symmetrization: Stability and a Sharp Energy
Threshold}
\maketitle
\begin{abstract}
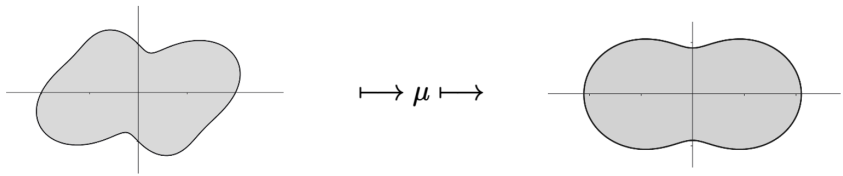
We develop a quantitative theory of the Brownian symmetrization for
simply connected planar domains with a prescribed real Brownian exit
law. We consider the coefficient energies
\[
\ensuremath{\mathcal{E}_{s}(U)=\sum_{n\ge1}n^{2s}|c_{n}|^{2}},
\]
 where 
\[
\sum_{n\ge1}c_{n}z^{n}
\]
 is a conformal map from the unit disk onto $U$ fixing the origin.
Building on the minimization principle contained in Pruss’s analytic
rearrangement theorem, we derive, for $0<s<1$, an explicit lower
bound for the energy deficit in terms of the difference of squared
conformal radii and prove strong convergence of normalized conformal
maps along minimizing sequences. A second main result identifies a
sharp threshold at the area exponent $s=\tfrac{1}{2}$. For $0<s<\tfrac{1}{2}$,
the ratio of the Gross energy to the original energy has a positive
lower bound depending only on $s$. For $\tfrac{1}{2}\le s\le1$,
every ratio in $(0,1]$ occurs at any prescribed positive energy,
already among bounded domains with analytic boundaries, with the exit
law allowed to vary. The proofs combine the exact periodized kernel
of the spectral fractional Laplacian, rearrangement inequalities,
and estimates exploiting univalence. We also establish sharp Wasserstein
convexity and support-diameter bounds for the minimum, characterize
finite minimal area through the logarithmic energy of the quantile’s
Stieltjes measure, and relate the uniform upper-energy index to moments
of the target law. In particular, we recover Popoviciu’s sharp variance
bound. We conclude by stating several open problems.
\end{abstract}

\smallskip{}

\section{Introduction}

Let $Z_{t}=X_{t}+iY_{t}$ be planar Brownian motion started at the
origin, with independent standard one-dimensional coordinates. For
a domain $U\subset\C$ containing $0$, write $\tau_{U}$ for its
first exit time. The planar Skorokhod embedding problem asks for a
domain such that $\Law(X_{\tau_{U}})$ equals a given centered probability
measure. Gross~\cite{Gross} constructed a canonical simply connected
solution from the quantile function of a nondegenerate law with finite
second moment. Boudabra and Markowsky~\cite{BMRemarks} extended
the construction to finite moments of order $p>1$ and established
an accompanying exit-time estimate. Their alternative construction~\cite{BMNew}
gives further solutions with different geometry. The endpoint $p=1$
under an additional Hilbert-transform integrability condition is treated
in~\cite{BoudabraEndpoint}.

The nonuniqueness of an embedding motivates optimization over domains
with the same real exit law. Becher, Boudabra, and Haggui~\cite[Theorem~5]{Energy}
established Gross minimization for the first-order coefficient functional
$\frac{1}{4}\sum n^{2}|c_{n}|^{2}$ of a conformal map $f(z)=\sum_{n\geq1}c_{n}z^{n}$
for which $U=f(\mathbb{D})$ is a solution for the planar Skorokhod
embedding problem. 

We consider the coefficient energies

\[
\mathcal{E}_{s}(U)=\sum_{n\ge1}n^{2s}|c_{n}|^{2},\qquad s\ge0,
\]
where $f(z)=\sum_{n\ge1}c_{n}z^{n}$ maps $\mathbb{D}$ conformally
onto $U$ and satisfies $f(0)=0$. These functionals are independent
of the choice of such a map, and $\mathcal{E}_{1/2}(U)=\operatorname{Area}(U)/\pi$.
They coincide with the squared Dirichlet norms considered by Pruss
\cite[§6]{Pruss}:
\[
\mathcal{E}_{s}(U)=\|f\|^{2}_{\mathcal{D}_{2s}}.
\]
We refer to $\mathcal{E}_{s}(U)$ as the boundary energy of order
$s$.

Our motivation is the planar Skorokhod embedding problem: we study
these energies over simply connected domains with a prescribed real
Brownian exit law. This formulation connects the analytic symmetrization
framework of Baernstein and Pruss \cite{Baernstein,Pruss} with the
probabilistic construction of Gross. The connection is explicit. If
$u$ is a real boundary datum with quantile $q$, its symmetric decreasing
rearrangement is $q(1-|\theta|/\pi)$, whereas the Gross construction
uses $q(|\theta|/\pi)$. The corresponding holomorphic extensions,
normalized to have the same value at the origin, differ by precomposition
with $z\mapsto-z$ and therefore have the same image \cite{Energy}.
In the Hardy-integrable setting, this identifies Gross symmetrization
with the analytic rearrangement appearing in Pruss’s theorem. Consequently,
the minimization inequality
\[
\mathcal{E}_{s}(U^{G}_{\mu})\le\mathcal{E}_{s}(U),\qquad0<s\le1,
\]
and its finite-energy equality characterization for $0<s<1$ are already
contained in \cite[Theorem 6.1]{Pruss}, after this identification.

We give a direct proof adapted to the exit-law formulation, using
the exact periodized kernel of the spectral fractional Laplacian \cite{CCM}.
Pruss’s proof integrates comparisons of quadratic forms associated
with resolvent convolution kernels. Our argument works directly with
the singular quadratic form, whose Fourier representation we establish
with explicit normalization constants, and applies the periodic rearrangement
approach recalled in \cite{StrictRearrangement}. Both arguments rely
on circular rearrangement. However, the direct kernel representation
provides the basis for our main results of this work about the quantitative
analysis: for $0<s<1$, we obtain an explicit estimate
\[
\mathcal{E}_{s}(U)-\mathcal{E}_{s}(U^{G}_{\mu})\ge\gamma_{s}\bigl(R^{2}_{\mu}-r(U)^{2}\bigr)\ge0,
\]
where $r(U)$ is the conformal radius at the origin, $R_{\mu}=r(U^{G}_{\mu})$,
and $\gamma_{s}>0$ depends only on $s$. Thus the energy deficit
controls the loss of conformal radius relative to the Gross domain.
We also establish compactness of energy sublevels in weaker coefficient
norms and convergence of minimizing sequences in the energy norm.
The treatment of equality is consistent with Pruss’s theorem and the
strict periodic rearrangement theory in \cite{StrictRearrangement},
with the signed boundary data required by the exit-law constraint
handled explicitly. Beyond minimization, we investigate the range
of the symmetrization ratio
\[
\rho_{s}(U)=\frac{\mathcal{E}_{s}(U^{G}_{\mu_{U}})}{\mathcal{E}_{s}(U)}.
\]
It turns out that its behavior changes at the area exponent $s=\tfrac{1}{2}$.
For $0<s<\tfrac{1}{2}$, univalence yields a positive lower bound
depending only on $s$. The proof combines weighted area estimates,
reconstruction from real boundary values, and a fractional Hardy estimate
controlling the relevant Lorentz integral. For $\tfrac{1}{2}\le s\le1$,
an explicit logarithmic family instead realizes every ratio in $(0,1]$,
after normalization to any prescribed positive energy, already among
bounded domains with analytic boundaries. In this range statement
the exit law varies with the domain. These results describe how much
energy symmetrization can remove, complementing the identification
of its minimizer.

The probabilistic formulation requires care with admissibility and
integrability. Nondegeneracy of the target law is necessary for the
Gross map to be univalent, while finite variance of the terminal real
coordinate alone does not imply finite mean exit time for a weak embedding.
We exhibit this obstruction and justify the required variance identities
through localization or Hardy integrability. The minimization inequality
nevertheless extends to the full simply connected admissible class:
finite competitor energy supplies the necessary $\mathbf{H}^{2}$
integrability, while infinite competitor energy makes the extended-value
inequality immediate.

We further study the dependence of the optimum on the target law.
Quantile interpolation yields an exact quadratic identity and a sharp
convexity bound in the quadratic Wasserstein distance, complementing
the Hardy and Wasserstein continuity results in \cite{BrownianConvergence}.
For bounded laws, we establish a sharp support-diameter bound when
$0<s<\tfrac{1}{2}$. We also express minimal area as the logarithmic
energy of the Stieltjes measure of the quantile, recovering the support-gap
obstruction discussed in \cite{BoundedGross} and characterizing the
existence of a finite-area simply connected embedding. Finally, the
subcritical estimates identify the uniform upper-energy index in terms
of the target law’s moment exponents and yield maximal-energy embeddings
below that index.

\section{Exit laws and admissible domains}

\label{pr:framework}

Let $Z_{t}=X_{t}+iY_{t}$ be planar Brownian motion started at $0$,
where $X$ and $Y$ are independent standard real Brownian motions.
Thus its generator is $\frac{1}{2}\Delta$. Throughout, a domain is
a connected open subset of $\C$, and every domain under consideration
is proper and contains $0$. Its first exit time and real exit law
are denoted by 
\[
\tau_{U}=\inf\{t\geq0:Z_{t}\notin U\},\qquad\mu_{U}=\operatorname{Law}(X_{\tau_{U}}),
\]
where the second expression is used when $\tau_{U}<\infty$ almost
surely. We write $\D=\{z\in\C:|z|<1\}$ and $\T=\partial\D$, and
use the probability measure $dm(\theta)=d\theta/(2\pi)$ on $\T$. 

The Hardy space $\mathbf{H}^{p}$ is the set of analytic functions
$f$ in the unit disc such that 
\[
\|f\|_{\mathbf{H}^{p}}:=\left(\sup_{0<r<1}\int^{\pi}_{-\pi}|f(re^{i\theta})|^{p}\,dm(\theta)\right)^{1/p}.
\]
The quantity $\|f\|_{\mathbf{H}^{p}}$ is called Hardy norm of $f$.
More details about the theory of Hardy spaces can be found in \cite{DurenHardy}.

\subsection{The admissible classes and the Gross construction}

Let $\mathcal{M}_{*}$ be the class of centered probability measures
$\mu$ on $\R$ such that $\mu\ne\delta_{0}$ and $\int_{\R}|x|^{p}\,d\mu(x)<\infty$
for some $p>1$. Let 
\[
\mathcal{M}_{2}=\left\{ \mu\in\mathcal{M}_{*}:\int_{\R}x^{2}\,d\mu(x)<\infty\right\} .
\]
For a centered law we allow $\Var(\mu)=\int_{\R}x^{2}\,d\mu(x)$ to
take the value $+\infty$. The nondegeneracy condition excludes the
constant analytic map arising from the quantile of $\delta_{0}$.
For $\mu\in\mathcal{M}_{*}$ define 
\begin{align*}
\mathcal{B}(\mu) & =\{U:\ \tau_{U}<\infty\text{ almost surely and }\mu_{U}=\mu\},\qquad\mathcal{B}=\bigcup_{\mu\in\mathcal{M}_{*}}\mathcal{B}(\mu)\\
\mathcal{C}(\mu) & =\{U\in\mathcal{B}(\mu):U\text{ is simply connected}\}.
\end{align*}
We call members of $\mathcal{B}(\mu)$ $\mu$-domains. No moment condition
on $\tau_{U}$ is imposed in this definition. When needed, a $\mu$-domain
is called strong if $\E\tau^{p/2}_{U}<\infty$ for some $p>1$. For
simply connected domains this is equivalent to the corresponding Riemann
map belonging to the Hardy space $\mathbf{H}^{p}$ for some $p>1$
\cite{Burkholder}. 

Fix $\mu\in\mathcal{M}_{*}$, write $F_{\mu}$ for its distribution
function, and define its quantile by 
\[
q_{\mu}(t)=\inf\{x\in\R:F_{\mu}(x)\geq t\},\qquad0<t<1.
\]
A key feature of a quantile is that it generates data sampling as
$\mu$ when fed with uniform points of $(0,1)$. The even periodic
function 
\[
\phi_{\mu}(\theta)=q_{\mu}(|\theta|/\pi),\qquad-\pi<\theta<\pi,
\]
has mean zero. Put 
\begin{equation}
a_{n}=\frac{1}{\pi}\int^{\pi}_{-\pi}\phi_{\mu}(\theta)\cos(n\theta)\,d\theta,\qquad G_{\mu}(z)=\sum_{n\geq1}a_{n}z^{n},\qquad U^{G}_{\mu}=G_{\mu}(\D).\label{pr:gross}
\end{equation}
Gross's construction~\cite{Gross} and its $p$-moment extension~\cite{BMRemarks}
assert that $G_{\mu}$ is univalent, $U^{G}_{\mu}\in\mathcal{C}(\mu)$,
and 
\begin{equation}
G_{\mu}\in\mathbf{H}^{p},\qquad\E\tau^{p/2}_{U^{G}_{\mu}}<\infty\quad\text{whenever}\quad1<p<\infty\ \text{ and }\ \int_{\R}|x|^{p}\,d\mu(x)<\infty.\label{pr:gross-moment}
\end{equation}
The boundary values satisfy 
\[
\operatorname{Re}G^{*}_{\mu}=\phi_{\mu}
\]
 almost everywhere \cite{Carleson}. The Brownian symmetrization on
$\mathcal{B}$ is defined by 
\begin{equation}
\mathfrak{B}(U)=U^{G}_{\mu_{U}}.\label{pr:brownian-symmetrization}
\end{equation}
Thus the definition uses the real exit law and applies to domains
of arbitrary connectivity within $\mathcal{B}$.

A first remark is that the first coefficient in \eqref{pr:gross}
has a universal property: a fixed sign independent of the choice of
$\mu$. Indeed, 
\begin{align*}
a_{1} & =2\int^{1}_{0}q_{\mu}(t)\cos(\pi t)\,dt\\
 & =\int^{1}_{0}\int^{1}_{0}\bigl(q_{\mu}(t)-q_{\mu}(v)\bigr)\bigl(\cos(\pi t)-\cos(\pi v)\bigr)\,dt\,dv<0.
\end{align*}
The integrand is nonpositive because the quantile is nondecreasing
and the cosine is strictly decreasing on $(0,1)$; strict inequality
follows because $\mu$ is nondegenerate. Consequently 
\begin{equation}
H_{\mu}(z):=G_{\mu}(-z),\qquad H_{\mu}'(0)=-a_{1}>0,\qquad H_{\mu}(\D)=U^{G}_{\mu}.\label{pr:gross-orientation}
\end{equation}
Its real boundary part is the symmetric decreasing rearrangement 
\[
\operatorname{Re}H^{*}_{\mu}(e^{i\theta})=q_{\mu}(1-|\theta|/\pi)\quad\text{for almost every }\theta\in(-\pi,\pi).
\]
Precomposition by $z\mapsto-z$ changes this orientation but does
not change the image domain.

\subsection{An analytic subclass and boundary laws}

Motivated by the covering-map framework used in symmetrization \cite{Pruss},
define $\mathcal{P}$ to be the subclass of $\mathcal{B}$ consisting
of domains admitting a holomorphic universal covering $f:\D\to U$
with $f(0)=0$ and $f\in\mathbf{H}^{p}$ for some $p>1$. Set 
\[
\mathcal{P}(\mu)=\mathcal{P}\cap\mathcal{B}(\mu).
\]
This is an explicit analytic definition; no assertion about a broader
classical symmetrization construction is needed. In the simply connected
case the covering is a Riemann map, and $\mathcal{C}(\mu)\cap\mathcal{P}$
is, in some sense, the strong subclass. The following lemma is standard
and can be found in \cite{DurenHardy}.

\begin{lemma}\label{pr:automorphisms} For $p>0$, composition with
an automorphism of $\D$ preserves membership in $\mathbf{H}^{p}$
in both directions. Consequently the Hardy condition in the definition
of $\mathcal{P}$ is independent of the chosen universal covering
and of its normalization at the base point. \end{lemma} 

We recall the boundary-law identification in the Hardy-integrable
setting; see also~\cite{BMRemarks,BoudabraEndpoint}. The proof specifies
the assumptions used below.

\begin{lemma}\label{pr:boundary-law} Let $U$ be a proper simply
connected domain containing $0$, with $\tau_{U}<\infty$ almost surely.
Let $f:\D\to U$ be a Riemann map with $f(0)=0$ and $f\in\mathbf{H}^{p}$
for some $p>1$. If $\Theta$ is uniform on $(-\pi,\pi)$, then 
\[
Z_{\tau_{U}}\ \stackrel{{\rm law}}{=}\ f^{*}(e^{i\Theta}),\qquad X_{\tau_{U}}\ \stackrel{{\rm law}}{=}\ \operatorname{Re}f^{*}(e^{i\Theta}).
\]
In particular, for $f(z)=\sum_{n\geq1}c_{n}z^{n}$ and $u=\operatorname{Re}f^{*}$,
\begin{equation}
\widehat{u}(0)=0,\qquad\widehat{u}(n)=\frac{c_{n}}{2},\qquad\widehat{u}(-n)=\frac{\overline{c_{n}}}{2}\quad(n\geq1),\qquad\widehat{u}(k)=\int^{\pi}_{-\pi}u(\theta)e^{-ik\theta}\,dm(\theta).\label{pr:boundary-fourier}
\end{equation}
\end{lemma} 

\begin{theorem}\label{pr:strict-inclusion} For every $\mu\in\mathcal{M}_{*}$
the class $\mathcal{P}(\mu)$ is nonempty. There exists $\mu\in\mathcal{M}_{*}$
for which $\mathcal{P}(\mu)\subsetneq\mathcal{B}(\mu)$. \end{theorem} 

Before delving into the proof, let's read what theorem \ref{pr:strict-inclusion}
says. The exit-law formulation extends the Hardy-class domain symmetrization
considered by Pruss \cite[§2]{Pruss}. On $\mathcal{P}$, the map
$U\mapsto U^{G}_{\mu_{U}}$ agrees with the Baernstein symmetrization
described there, while it remains well defined on the larger class
$\mathcal{B}$, without requiring Hardy integrability of a universal
covering map. The strict inclusion $\mathcal{P}(\mu)\subsetneq\mathcal{B}(\mu)$
for some $\mu\in\mathcal{M}_{*}$, even with compact support, shows
that this enlargement is substantive. Thus Brownian symmetrization
applies to domains outside that Hardy-class framework, while its output
belongs to $\mathcal{P}(\mu)$ and preserves the prescribed real exit
law.
\begin{proof}
Nonemptiness follows from \eqref{pr:gross-moment}. For strictness,
consider 
\[
U=\C\setminus([-2,2]-i),\qquad J(\zeta)=\zeta+\zeta^{-1}-i.
\]
The map $J$ is a conformal bijection from $\{|\zeta|>1\}$ onto $U$.
Put $a=(1+\sqrt{5})/2$, so that $J(ia)=0$. Brownian motion in the
exterior disk, started at $ia$, reaches its boundary almost surely.
Under $J$ the conformal time change is finite up to this hitting
time, since $|J'(\zeta)|\leq2$ for $|\zeta|>1$. It follows that
Brownian motion in $U$, started at $0$, exits almost surely. Its
real exit law $\mu$ is supported on $[-2,2]$. Reflection in the
imaginary axis preserves $U$ and the starting point, so $\mu$ is
centered. The exterior-disk Poisson kernel is positive on the circle,
and its image under $J$ gives positive mass to both halves of the
slit. Thus $\mu$ is nondegenerate and belongs to $\mathcal{M}_{*}$.

Let 
\[
h(z)=\frac{1+z}{1-z},\qquad F(z)=J(e^{h(z)})=2\cosh(h(z))-i.
\]
The map $h$ sends $\D$ conformally onto the right half-plane, and
the exponential is a universal covering of $\{|\zeta|>1\}$. Hence
$F$ is a universal covering of $U$. For $|\theta|\leq1-r$, 
\[
|1-re^{i\theta}|\leq2(1-r),\qquad\operatorname{Re}h(re^{i\theta})=\frac{1-r^{2}}{|1-re^{i\theta}|^{2}}\geq\frac{1}{4(1-r)}.
\]
For $r$ sufficiently close to $1$, $|e^{h(re^{i\theta})}|\geq3$
on this arc. Since $|J(\zeta)|\geq|\zeta|/2$ when $|\zeta|\geq3$,
for each $p>0$ we obtain 
\[
\int^{\pi}_{-\pi}|F(re^{i\theta})|^{p}\,d\theta\geq2^{1-p}(1-r)\exp\!\left(\frac{p}{4(1-r)}\right)\longrightarrow\infty.
\]
Thus $F\notin\mathbf{H}^{p}$ for every $p>0$.

To make the base point explicit, set 
\[
w_{0}=\log a+\frac{\pi i}{2},\qquad z_{0}=\frac{w_{0}-1}{w_{0}+1},\qquad\alpha(z)=\frac{z+z_{0}}{1+\overline{z_{0}}z}.
\]
We have $z_{0}\in\D$, $e^{h(z_{0})}=ia$, and $(F\circ\alpha)(0)=J(ia)=0$.
Lemma~\ref{pr:automorphisms} shows that this based covering also
belongs to no $\mathbf{H}^{p}$, and the same is true of every other
universal covering of $U$. Therefore $U\in\mathcal{B}(\mu)\setminus\mathcal{P}(\mu)$. 
\end{proof}

\begin{remark}\label{pr:slit-law} The law of $X_{\tau_{U}}$in Theorem~\ref{pr:strict-inclusion}
can be written explicitly: 
\[
d\mu(x)=\frac{\sqrt{5}}{\pi(1+x^{2})\sqrt{4-x^{2}}}\mathbf{1}_{(-2,2)}(x)\,dx.
\]
Indeed, writing $r=1/a$, exterior-disk harmonic measure at $ia$
has angle density $(1-r^{2})/[2\pi(1+r^{2}-2r\sin\theta)]$. Pushforward
by $x=2\cos\theta$, summing the two preimages and using $1-r^{2}=r$,
gives the displayed density. Consequently 
\[
F_{\mu}(x)=\frac{1}{2}+\frac{1}{\pi}\arctan\!\left(\frac{\sqrt{5}x}{\sqrt{4-x^{2}}}\right)\quad(-2<x<2),
\]
and 
\[
q_{\mu}(t)=\frac{2\sin(\pi(t-\tfrac{1}{2}))}{\sqrt{1+4\cos^{2}(\pi(t-\tfrac{1}{2}))}},\qquad0<t<1.
\]

\end{remark}

\subsection{Exit-time identities and their hypotheses}

The coefficient identity below is classical and is recalled by Gross~\cite{Gross},
with attribution to Bañuelos and Carroll \cite{Ba=0000F1uelos}. We
include the exhaustion proof to specify its extended-value interpretation
and the hypotheses for the variance identity.

\begin{lemma}\label{pr:exit-time} Let $U$ be a proper simply connected
domain containing $0$, and let $f(z)=\sum_{n\geq1}c_{n}z^{n}$ map
$\D$ conformally onto $U$, with $f(0)=0$. Then, in $[0,\infty]$,
\begin{equation}
\E\tau_{U}=\frac{1}{2}\sum_{n\geq1}|c_{n}|^{2}.\label{pr:exit-coefficients}
\end{equation}
If this quantity is finite, then $\tau_{U}<\infty$ almost surely
and, writing $\mu=\mu_{U}$, one has 
\begin{equation}
\E\tau_{U}=\Var(\mu).\label{pr:exit-variance}
\end{equation}
More generally, \eqref{pr:exit-variance} holds in the extended nonnegative
sense for $U\in\mathcal{C}(\mu)\cap\mathcal{P}$. \end{lemma} 
\begin{proof}
Use the exhaustion $U_{r}=f(r\D)$ and its exit times $\tau_{r}$
from Lemma~\ref{pr:boundary-law}; the exhaustion itself does not
require a Hardy hypothesis. Each $U_{r}$ is bounded, so comparison
with a containing disk gives $\E\tau_{r}<\infty$. Apply optional
stopping first at $t\wedge\tau_{r}$ and then let $t\to\infty$; the
stopped spatial process is bounded, and $\E(t\wedge\tau_{r})\uparrow\E\tau_{r}$.
Applying this to $|Z_{t}|^{2}-2t$, and then using conformal invariance
of harmonic measure and Parseval's identity on the inner circle, gives
\[
\E\tau_{r}=\frac{1}{2}\E|Z_{\tau_{r}}|^{2}=\frac{1}{2}\int^{\pi}_{-\pi}|f(re^{i\theta})|^{2}\,dm(\theta)=\frac{1}{2}\sum_{n\geq1}|c_{n}|^{2}r^{2n}.
\]
Since $\tau_{r}\uparrow\tau_{U}$, monotone convergence on both sides
proves \eqref{pr:exit-coefficients}.

Suppose now that $\E\tau_{U}<\infty$. The stopped stochastic integrals
$X_{t\wedge\tau_{U}}=\int^{t}_{0}\mathbf{1}_{\{s<\tau_{U}\}}\,dX_{s}$
converge in $L^{2}$ to $X_{\tau_{U}}$. The Itô isometry~\cite[Chapter~7]{MortersPeres}
gives 
\[
\E X_{\tau_{U}}=0,\qquad\E X^{2}_{\tau_{U}}=\E\int^{\infty}_{0}\mathbf{1}_{\{s<\tau_{U}\}}\,ds=\E\tau_{U},
\]
which proves \eqref{pr:exit-variance} under the finite-mean hypothesis.

Finally suppose $f\in\mathbf{H}^{p}$ for some $p>1$, as holds for
$U\in\mathcal{C}(\mu)\cap\mathcal{P}$. Set $u=\operatorname{Re}f^{*}$.
If $\Var(\mu)<\infty$, Lemma~\ref{pr:boundary-law} gives $u\in L^{2}$
and \eqref{pr:boundary-fourier}. Parseval's identity therefore yields
\[
\Var(\mu)=\int^{\pi}_{-\pi}|u|^{2}\,dm=\frac{1}{2}\sum_{n\geq1}|c_{n}|^{2}=\E\tau_{U}.
\]
If $\Var(\mu)=\infty$, the finite-mean assertion already proved rules
out $\E\tau_{U}<\infty$. Thus both sides are infinite. 
\end{proof}

\begin{corollary}\label{pr:strong-distance} For $U\in\mathcal{C}(\mu)\cap\mathcal{P}$,
\[
\frac{\operatorname{dist}(0,\partial U)^{2}}{2}\leq\E\tau_{U}=\Var(\mu),
\]
with the usual interpretation when the variance is infinite. \end{corollary} 
\begin{proof}
The disk centered at $0$ with radius $\operatorname{dist}(0,\partial U)$
is contained in $U$ and has mean exit time equal to half the square
of that radius. Domain inclusion and Lemma~\ref{pr:exit-time} give
the conclusion. 
\end{proof}

The necessity of an exit-moment condition is already illustrated in~\cite{BMRemarks}.
The following variation fixes the exit law and the distance to the
boundary explicitly.

\begin{proposition}\label{pr:weak-counterexample} A finite second
moment of $\mu$ alone does not imply \eqref{pr:exit-variance}, even
for simply connected $\mu$-domains. More precisely, for $a,b>0$
let 
\[
V_{a,b}=\C\setminus\{x+iy:x\in\{-a,a\},\ |y|\geq b\}.
\]
Then 
\[
V_{a,b}\in\mathcal{C}\left(\tfrac{1}{2}(\delta_{-a}+\delta_{a})\right),\qquad\E\tau^{p/2}_{V_{a,b}}=\infty\quad(p\geq1),\qquad\operatorname{dist}(0,\partial V_{a,b})^{2}=a^{2}+b^{2}.
\]
\end{proposition} 
\begin{proof}
The horizontal corridor $|\operatorname{Im}z|<b$ connects all components
between the removed rays, so $V_{a,b}$ is connected. Its complement
in the Riemann sphere consists of four rays joined at infinity and
is connected; hence $V_{a,b}$ is simply connected. Planar Brownian
motion hits any fixed nontrivial closed subsegment of a removed ray
almost surely, by recurrence and the positive logarithmic capacity
of a line segment~\cite[Chapters~3 and~8]{MortersPeres}. Thus $\tau_{V_{a,b}}$
is finite almost surely. Its real exit coordinate is $-a$ or $a$;
reflection in the imaginary axis gives probability $1/2$ for each.

Let $\sigma=\inf\{t\geq0:X_{t}=2a\}$. With positive probability the
path remains in $|Y_{t}|<b$ up to $\sigma$, and on this event $\sigma<\tau_{V_{a,b}}$.
Positivity follows, for example, from the independence of $X$ and
$Y$ by intersecting $\{\sigma<T\}$ with $\{\sup_{t\leq T}|Y_{t}|<b\}$
for any fixed $T>0$. Set $E=\{\sigma<\tau_{V_{a,b}}\}\in\mathcal{F}_{\sigma}$;
the preceding argument gives $\mathbb{P}(E)>0$. On $E$, after $\sigma$,
the path cannot leave $V_{a,b}$ before its real coordinate returns
to $a$, because the half-plane $\{\operatorname{Re}z>a\}$ is contained
in $V_{a,b}$. Conditional on $\mathcal{F}_{\sigma}$, the strong
Markov property gives this return time the law of $H_{a}$, the time
for one-dimensional Brownian motion to move a distance $a$ in a prescribed
direction. The reflection principle~\cite[Chapter~2]{MortersPeres}
implies $\mathbb{P}(H_{a}>t)\sim a\sqrt{2/\pi}\,t^{-1/2}$ as $t\to\infty$.
Consequently, by conditioning on $\mathcal{F}_{\sigma}$ on $E$ and
using nonnegative expectations, 
\[
\E\tau^{q}_{V_{a,b}}\geq\mathbb{P}(E)\,\E H^{q}_{a}=\infty,\qquad q\geq1/2.
\]
Finally the closest boundary points are the four ray endpoints, each
at distance $\sqrt{a^{2}+b^{2}}$ from $0$. 
\end{proof}

\begin{remark}\label{pr:weak-consequences} For the domains in Proposition~\ref{pr:weak-counterexample},
the variance of the real exit law is $a^{2}$, although the mean exit
time is infinite. If $b>a$, even the bound $\operatorname{dist}(0,\partial V_{a,b})^{2}/2\leq\Var(\mu_{V_{a,b}})$
fails. These domains are also symmetric with respect to the real axis
and have connected vertical sections. Thus such geometric properties
alone cannot replace the exit-moment hypothesis in the uniqueness
theorem for the Gross domain \cite{BMRemarks}. \end{remark}

\section{Spectral boundary energies and rearrangement}

\label{sp:section}

For a proper simply connected domain $U$ containing $0$, choose
a conformal map $f_{U}:\D\to U$ with $f_{U}(0)=0$ and write $f_{U}(z)=\sum_{n\ge1}c_{n}z^{n}$.
For $s\ge0$, define the fractional boundary energy of $U$ (of index
$s$) by
\begin{equation}
\mathcal{E}_{s}(U):=\sum_{n\ge1}n^{2s}|c_{n}|^{2}\in(0,\infty].\label{sp:energy-definition}
\end{equation}
Any two such maps differ by a rotation of the disk, so their coefficients
have the same absolute values. Thus this is a well-defined functional
of the domain with its specified starting point. The area formula
for univalent maps~\cite{DurenUnivalent} gives 
\[
A(U):=\operatorname{Area}(U)=\pi\mathcal{E}_{1/2}(U);
\]
 the energy studied in~\cite{Energy} is $\frac{1}{4}\mathcal{E}_{1}(U)$.
These coefficient functionals are used only on simply connected domains.
Translating a domain while keeping the starting point at $0$ generally
changes them; area is an exception.

We identify the unit circle with $\R/(2\pi\mathbb{Z})$ and use angular
Lebesgue measure in all integrals. For $u\in L^{2}(-\pi,\pi)$, extended
periodically, set 
\begin{equation}
\widehat{u}_{k}=\frac{1}{2\pi}\int^{\pi}_{-\pi}u(\theta)e^{-ik\theta}\,d\theta,\qquad\mathscr{S}_{s}(u)=\sum_{k\ne0}|k|^{2s}|\widehat{u}_{k}|^{2}.\label{sp:fourier}
\end{equation}
The sum is allowed to be infinite. We define the periodic Sobolev
space $H^{s}(\mathbb{T})$ by

\[
H^{s}(\mathbb{T}):=\{u\in L^{2}(\mathbb{T}):\mathscr{S}_{s}(u)<\infty\}.
\]
We equip this space with the norm
\[
\|u\|_{H^{s}(\mathbb{T})}:=\left(\|u\|^{2}_{L^{2}(\mathbb{T})}+\mathscr{S}_{s}(u)\right)^{1/2}.
\]

The spectral fractional Laplacian on the torus has an exact periodized
power kernel; see Roncal and Stinga~\cite{RoncalStinga}. We record
its one-dimensional quadratic form with the Fourier normalization
used here. For $0<s<1$, define 
\begin{equation}
K_{s}(t)=\sum_{j\in\mathbb{Z}}|t+2\pi j|^{-1-2s},\qquad c_{s}=\frac{\Gamma(1+2s)\sin(\pi s)}{\pi},\qquad t\notin2\pi\mathbb{Z},\label{sp:kernel}
\end{equation}
and 
\begin{equation}
Q_{s}(u)=\int^{\pi}_{-\pi}\int^{\pi}_{-\pi}|u(\theta)-u(\sigma)|^{2}K_{s}(\theta-\sigma)\,d\sigma\,d\theta\label{sp:quadratic}
\end{equation}
the underlying quadratic form. The integrand is assigned the value
zero on the diagonal, a set of two-dimensional measure zero. The next
lemma is elementary.

\begin{lemma}\label{sp:monotonicity} The kernel $K_{s}$ is even,
$2\pi$-periodic, and strictly decreasing on $(0,\pi)$. \end{lemma} 

It turns out that the quantity $\mathscr{S}_{s}$ is proportional
to $Q_{s}$. More precisely

\begin{proposition} \cite{Ohrem} \label{sp:identity} For $0<s<1$
and $u\in L^{2}(\T)$, 
\begin{equation}
\mathscr{S}_{s}(u)=\frac{c_{s}}{4\pi}Q_{s}(u),\label{sp:exact}
\end{equation}
with equality in the extended nonnegative sense. \end{proposition} 

We use the circular Riesz rearrangement inequality: if $f,h$ are
nonnegative measurable periodic functions and $J$ is nonnegative,
even, periodic, and nonincreasing on $(0,\pi)$, then 
\begin{equation}
\iint f(\theta)J(\theta-\sigma)h(\sigma)\,d\theta\,d\sigma\leq\iint f^{\#}(\theta)J(\theta-\sigma)h^{\#}(\sigma)\,d\theta\,d\sigma.\label{sp:riesz}
\end{equation}
Both double integrals are over $(-\pi,\pi)^{2}$. This is the circle
inequality discussed in the published paper~\cite{CCM}.

\begin{theorem}\label{sp:rearrangement} For every real $u\in L^{2}(\T)$
and $0<s<1$, 
\begin{equation}
\mathscr{S}_{s}(u^{\#})\leq\mathscr{S}_{s}(u).\label{sp:contraction}
\end{equation}
\end{theorem} 

As explained in the introduction, Theorem \ref{sp:rearrangement}
is a reformulation of Pruss’s contraction theorem \cite{Pruss}. We
include an adapted proof using the explicit periodized kernel and
the periodic rearrangement approach of Cabré--Csató--Mas \cite{CCM}.
Shifts and truncations allow us to apply the nonnegative rearrangement
inequality to signed boundary data, obtaining an equimeasurable rearrangement
that preserves the prescribed real boundary distribution. The exact
Fourier identity then converts the quadratic-form comparison into
the corresponding spectral-energy inequality. We retain this argument
because it applies also to the modified kernels used below to derive
the quantitative deficit estimate.
\begin{proof}
First suppose that $v$ is bounded. For $N>0$ let $J_{N}=\min\{K_{s},N\}$.
By Lemma~\ref{sp:monotonicity}, this is a bounded nonnegative kernel
satisfying the hypotheses of \eqref{sp:riesz}. Choose $M\geq\|v\|_{\infty}$
and apply \eqref{sp:riesz} with $f=h=v+M$. Since $(v+M)^{\#}=v^{\#}+M$,
expansion of the quadratic form gives 
\[
\iint|v^{\#}(\theta)-v^{\#}(\sigma)|^{2}J_{N}(\theta-\sigma)\,d\theta\,d\sigma\leq\iint|v(\theta)-v(\sigma)|^{2}J_{N}(\theta-\sigma)\,d\theta\,d\sigma.
\]
Indeed, for any integrable periodic kernel $J$ the expansion is 
\[
2\left(\int^{\pi}_{-\pi}J\right)\left(\int^{\pi}_{-\pi}v^{2}\right)-2\iint v(\theta)J(\theta-\sigma)v(\sigma)\,d\theta\,d\sigma.
\]
Equimeasurability preserves the first term and every term introduced
by the constant shift. Letting $N\to\infty$ by monotone convergence
proves $Q_{s}(v^{\#})\leq Q_{s}(v)$.

For general real $u\in L^{2}$, put $T_{M}(t)=\max\{-M,\min\{t,M\}\}$
and $u_{M}=T_{M}(u)$. The superlevel-set definition implies $(u_{M})^{\#}=T_{M}(u^{\#})$
almost everywhere. For fixed real $a,b$, $|T_{M}(a)-T_{M}(b)|$ increases
to $|a-b|$ as $M\to\infty$. Consequently monotone convergence on
both sides of $Q_{s}((u_{M})^{\#})\leq Q_{s}(u_{M})$ yields $Q_{s}(u^{\#})\leq Q_{s}(u)$.
Lemma~\ref{sp:identity} now proves \eqref{sp:contraction}. 
\end{proof}

\begin{corollary}\label{sp:domain} Let $\mu\ne\delta_{0}$ be a
centered probability measure with a finite $p$th moment for some
$p>1$, and let $U^{G}_{\mu}=G_{\mu}(\D)$ be its Gross domain. If
$U$ is a proper simply connected $\mu$-domain containing $0$, then
\begin{equation}
\mathcal{E}_{s}(U^{G}_{\mu})\leq\mathcal{E}_{s}(U),\qquad0<s\leq1,\label{sp:domainineq}
\end{equation}
where both energies are allowed to be infinite. In particular, the
Gross domain attains the infimum over this simply connected class,
and $\operatorname{Area}(U^{G}_{\mu})\leq\operatorname{Area}(U)$.
\end{corollary}

\begin{remark}\label{sp:infinitevariance} The case $s=1$ follows
by letting $s\uparrow1$, which represents a different proof than
in \cite{Energy}. The first-order minimization was established in~\cite{Energy}.
If $\mu$ has infinite second moment, the finite-energy reduction
in the proof shows that $\mathcal{E}_{s}(U)=\infty$ for every competitor
and every $s>0$. The extended minimization statement remains valid,
but has no finite-energy competitor in that case. \end{remark}

\section{Rigidity and an explicit energy deficit}

\label{rg:section}

For a simply connected domain $U$ containing $0$, its conformal
radius at $0$ is $r(U)=f_{U}'(0)$, where the Riemann map is normalized
by $f_{U}(0)=0$ and $f_{U}'(0)>0$. Write 
\[
R_{\mu}=r(U^{G}_{\mu})=H_{\mu}'(0)=-2\int^{1}_{0}q_{\mu}(t)\cos(\pi t)\,dt>0,
\]
where $H_{\mu}$ is defined in \eqref{pr:gross-orientation}. We shall
call $R_{\mu}$ the \emph{Gross conformal radius} of $\mu$. The equality
cases for circular Riesz rearrangement and periodic Gagliardo seminorms
are treated in~\cite{CCM} and~\cite{StrictRearrangement}, respectively.
The latter uses rearrangement of the absolute value. We give a direct
argument for the signed data used here; it also gives an explicit
domain deficit.

The next lemma is a signed form of the Hardy--Littlewood rearrangement
inequality, with $u^{\#}$ denoting the signed equimeasurable rearrangement
of $u$. Since both $u$ and the cosine weight may change sign, the
usual formulation for nonnegative functions is not directly applicable.
We provide a direct proof using the signed layer-cake representation,
the boundedness of the cosine weight, and its strict decrease on $(0,\pi)$.
This argument applies to every real $u\in L^{1}(\mathbb{T})$ without
truncation and also establishes the equality characterization.

\begin{lemma}\label{rg:linear} Let $u\in L^{1}(\T;\R)$ and set
$v=u^{\#}$. Then 
\begin{equation}
\int^{\pi}_{-\pi}u(\theta)\cos\theta\,d\theta\leq\int^{\pi}_{-\pi}v\cos\theta\,d\theta.\label{rg:linear-bound}
\end{equation}
Equality holds if and only if $u=v$ almost everywhere. \end{lemma} 
\begin{proof}
Among sets of a fixed measure $2a\in(0,2\pi)$, the centered arc $I=(-a,a)$
uniquely maximizes the integral of $\cos\theta$. Indeed, for any
such set $A$, 
\[
\int_{I}\cos\theta\,d\theta-\int_{A}\cos\theta\,d\theta=\int^{\pi}_{-\pi}(\mathbf{1}_{I}-\mathbf{1}_{A})(\cos\theta-\cos a)\,d\theta\geq0.
\]
The integrand is positive on $I\mathbin{\triangle}A=I\cup A\setminus I\cap A$,
apart from its null boundary, so equality forces $A=I$ modulo null
sets. The cases of measure $0$ and $2\pi$ are immediate. Apply this
to $A_{t}=\{u>t\}$ and $I_{t}=\{v>t\}$ for every real $t$. The
signed layer-cake identity gives 
\[
\int^{\pi}_{-\pi}(v-u)\cos\theta\,d\theta=\int_{\R}\left(\int_{I_{t}}\cos\theta\,d\theta-\int_{A_{t}}\cos\theta\,d\theta\right)dt.
\]
Fubini is justified by $\int_{\R}\int_{\T}|\mathbf{1}_{\{u>t\}}-\mathbf{1}_{\{v>t\}}|\,d\theta\,dt=\|u-v\|_{L^{1}(d\theta)}<\infty$.
The integrand on the right is nonnegative. If its integral is zero,
then $A_{t}=I_{t}$ modulo null sets for almost every $t$; the same
identity for the absolute difference implies $u=v$ almost everywhere. 
\end{proof}

\begin{theorem}\label{rg:radius} Let $\mu\in\mathcal{M}_{*}$ and
$U\in\mathcal{C}(\mu)\cap\mathcal{P}$. Then 
\begin{equation}
r(U)\leq R_{\mu},\label{rg:radius-bound}
\end{equation}
with equality if and only if $U=U^{G}_{\mu}$. In other words, Gross
domain maximizes the conformal radius at the origin among admissible
domains with the same real exit law $\mu$.\end{theorem} 
\begin{proof}
The normalized map $f_{U}$ belongs to $\mathbf{H}^{p}$ for some
$p>1$. By Lemma~\ref{pr:boundary-law}, $u=\operatorname{Re}f^{*}_{U}$
has law $\mu$, and its rearrangement is $v=\operatorname{Re}H^{*}_{\mu}$.
The Fourier identities give 
\[
\int^{\pi}_{-\pi}u\cos\theta\,d\theta=\pi r(U),\qquad\int^{\pi}_{-\pi}v\cos\theta\,d\theta=\pi R_{\mu}.
\]
Lemma~\ref{rg:linear} proves the inequality. In the equality case,
$u=v$ almost everywhere. Their boundary Fourier coefficients and
the normalizations at $0$ then give $f_{U}=H_{\mu}$, hence $U=U^{G}_{\mu}$.
The converse is immediate. 
\end{proof}

\begin{remark} During the preparation of this work, we were thinking
that the energy comparison does persist on concentric subdisks. But
in virtue of theorem \ref{rg:radius}, this thought is wrong. Indeed,
let $f_{U}(z)=\sum_{n\ge1}b_{n}z^{n}$ and $G_{\mu}(z)=\sum_{n\ge1}a_{n}z^{n}$,
and suppose that $U$ satisfies the hypotheses of the conformal-radius
theorem and $U\neq U^{G}_{\mu}$. Then

\[
\sum_{n\ge1}n^{2s}\bigl(|a_{n}|^{2}-|b_{n}|^{2}\bigr)r^{2n}=\bigl(R^{2}_{\mu}-r(U)^{2}\bigr)r^{2}+O(r^{4})>0
\]

for all sufficiently small $r>0$. Thus the global contraction coexists
with a strict increase of the energy on sufficiently small concentric
subdisks. This comparison concerns radial restrictions of the original
maps; rearranging the boundary values of $f_{U}(r\,\cdot)$ produces
a generally different map.

\end{remark}

We next establish a quantitative refinement of the minimization inequality.
The following theorem bounds the energy deficit from below by an explicit
multiple of the difference between the squared conformal radii of
the Gross domain and the competing domain.

\begin{theorem}\label{rg:deficit} Fix $0<s<1$, put $\alpha=1+2s$,
and define 
\begin{equation}
\varepsilon_{s}=2\alpha(\alpha+1)\pi^{-\alpha-2}(1-2^{-\alpha-2})\zeta(\alpha+2),\qquad\gamma_{s}=\pi c_{s}\varepsilon_{s}>0\label{rg:constant}
\end{equation}
with $c_{s}=\frac{\Gamma(1+2s)\sin(\pi s)}{\pi}$. If $U\in\mathcal{C}(\mu)$
and $\mathcal{E}_{s}(U)<\infty$, then 
\begin{equation}
\mathcal{E}_{s}(U)-\mathcal{E}_{s}(U^{G}_{\mu})\geq\gamma_{s}\bigl(R^{2}_{\mu}-r(U)^{2}\bigr)\geq0.\label{rg:energy-deficit}
\end{equation}
Consequently equality in the energy minimization holds if and only
if $U=U^{G}_{\mu}$. In particular, if $\mathcal{E}_{s}(U^{G}_{\mu})<\infty$,
then $U^{G}_{\mu}$ is the unique minimizer on $\mathcal{C}(\mu)$.
For finite-area competitors, 
\begin{equation}
A(U)-A(U^{G}_{\mu})\geq\frac{\pi}{8}\bigl(R^{2}_{\mu}-r(U)^{2}\bigr).\label{rg:area-deficit}
\end{equation}
\end{theorem} 
\begin{proof}
Twice differentiating the paired series for $K_{s}$ gives 
\[
K_{s}''(t)=\alpha(\alpha+1)\sum_{j\geq0}\left((2\pi j+t)^{-\alpha-2}+(2\pi(j+1)-t)^{-\alpha-2}\right).
\]
Termwise differentiation is justified by local uniform convergence
away from $2\pi\mathbb{Z}$. Convexity of $x\mapsto x^{-\alpha-2}$,
applied to each pair with midpoint $(2j+1)\pi$, proves $K_{s}''(t)\geq K_{s}''(\pi)=\varepsilon_{s}$.
Symmetry gives $K_{s}'(\pi)=0$, and therefore, for $0<t<\pi$, 
\[
-K_{s}'(t)=\int^{\pi}_{t}K_{s}''(x)\,dx\geq\varepsilon_{s}(\pi-t)\geq\varepsilon_{s}\sin t.
\]
It follows that 
\[
J_{s}(t)=\varepsilon_{s}(1+\cos t),\qquad\widetilde{K}_{s}(t)=K_{s}(t)-J_{s}(t)
\]
are nonnegative, even, periodic, and nonincreasing on $(0,\pi)$:
for $\widetilde{K}_{s}$ the derivative inequality proves monotonicity,
while $\widetilde{K}_{s}(\pi)=K_{s}(\pi)>0$ proves nonnegativity.
Finite energy implies $f_{U}\in\mathbf{H}^{2}$. Set $u=\operatorname{Re}f^{*}_{U}$
and $v=u^{\#}=\operatorname{Re}H^{*}_{\mu}$. For a nonnegative kernel
$J$, write 
\[
Q_{J}(w)=\iint|w(\theta)-w(\sigma)|^{2}J(\theta-\sigma)\,d\theta\,d\sigma.
\]
The truncation and constant-shift proof of Proposition~\ref{sp:rearrangement}
applies to $J_{s}$ and $\widetilde{K}_{s}$, since it uses only their
nonnegativity, evenness, and monotonicity. Thus $Q_{\widetilde{K}_{s}}(u)\geq Q_{\widetilde{K}_{s}}(v)$
and 
\[
Q_{s}(u)-Q_{s}(v)\geq Q_{J_{s}}(u)-Q_{J_{s}}(v).
\]
All these quadratic forms are finite: $J_{s},\widetilde{K}_{s}\leq K_{s}$
and the contraction theorem bounds the energy of $v$. Expanding the
bounded kernel form, and using the equal means and squared norms of
$u,v$, gives 
\begin{align*}
Q_{J_{s}}(u)-Q_{J_{s}}(v) & =2\varepsilon_{s}\left(\left|\int^{\pi}_{-\pi}v(\theta)e^{i\theta}\,d\theta\right|^{2}-\left|\int^{\pi}_{-\pi}u(\theta)e^{i\theta}\,d\theta\right|^{2}\right)\\
 & =2\varepsilon_{s}\pi^{2}\bigl(R^{2}_{\mu}-r(U)^{2}\bigr).
\end{align*}
The exact identity \eqref{sp:exact} and the relation $\mathcal{E}_{s}(U)=2\mathscr{S}_{s}(u)$
now yield \eqref{rg:energy-deficit}. Theorem~\ref{rg:radius} proves
its final inequality and the equality characterization. If the Gross
energy is finite, an infinite-energy competitor cannot be a minimizer.
Finally, at $s=1/2$ one has $c_{s}=1/\pi$, $\alpha=2$, and 
\[
\varepsilon_{1/2}=12\pi^{-4}(1-2^{-4})\zeta(4)=\frac{1}{8}.
\]
Multiplication by $\pi$ proves \eqref{rg:area-deficit}. 
\end{proof}

\begin{remark}\label{rg:scope} The finite-energy qualification excludes
the uninformative equality $\infty=\infty$. The constant in \eqref{rg:energy-deficit}
is explicit; we do not assert it is optimal for domain competitors.
This argument classifies equality only for $0<s<1$; the endpoint
$s=1$ contraction does not itself give an equality classification.
\end{remark}

\section{Stability and variation of the prescribed law}

\label{tr:section}

We first study compactness for a fixed exit law, then allow that law
to vary. For an analytic function $h(z)=\sum_{n\geq1}b_{n}z^{n}$,
recall the Dirichlet energy of $h$
\[
\|h\|^{2}_{\mathcal{D}_{2s}}:=\sum_{n\geq1}n^{2s}|b_{n}|^{2},\qquad s\geq0
\]
and 
\[
\mathcal{D}_{2s}=\left\{ h|\|h\|^{2}_{\mathcal{D}_{2s}}<\infty,h(0)=0\right\} ,\qquad s\geq0
\]
Thus $\mathcal{D}_{0}=\mathbf{H}^{2}\cap\{h:h(0)=0\}$, and for a
Riemann map $f_{U}$ fixing $0$ one has $\|f_{U}\|^{2}_{\mathcal{D}_{2s}}=\mathcal{E}_{s}(U)$.
Whenever maps of different domains are compared, we impose the normalization
$f_{U}'(0)>0$.

\subsection{Compactness and convergence of minimizing sequences}

\begin{proposition}\label{tr:compactness} Let $\mu\in\mathcal{M}_{2}$,
$s>0$, and $M<\infty$. The set 
\[
\mathcal{D}_{2s,M}(\mu)=\{f_{U}:U\in\mathcal{C}(\mu),\ f_{U}(0)=0,\ f_{U}'(0)>0,\ \mathcal{E}_{s}(U)\leq M\}
\]
is compact in $\mathcal{D}_{2t}$ for every $0\leq t<s$. If a sequence
in this set converges in $\mathbf{H}^{2}$ to $f$, then $f=f_{U}$
for a domain $U\in\mathcal{C}(\mu)$ and 
\[
\mathcal{E}_{s}(U)\leq\liminf_{j\to\infty}\mathcal{E}_{s}(U_{j}).
\]
\end{proposition} 
\begin{proof}
Write $f_{j}(z)=\sum_{n\geq1}c_{j,n}z^{n}$. The energy bound gives
\begin{equation}
\sum_{n>N}n^{2t}|c_{j,n}|^{2}\leq(N+1)^{-2(s-t)}M,\qquad0\leq t<s.\label{tr:tail}
\end{equation}
The coefficients are bounded, so diagonal selection gives a subsequence
converging coefficientwise to $(c_{n})$. Fatou's lemma gives $\sum n^{2s}|c_{n}|^{2}\leq M$.
The tail estimate and convergence of every finite coefficient sum
then imply convergence to $f(z)=\sum c_{n}z^{n}$ in $\mathcal{D}_{2t}$
for every $t<s$, in particular in $\mathbf{H}^{2}$ and locally uniformly
in $\D$. Every $f_{j}$ belongs to $\mathbf{H}^{2}$, and the exit-law
identity gives 
\[
\|f_{j}\|^{2}_{\mathbf{H}^{2}}=2\Var(\mu)>0.
\]
Strong $\mathbf{H}^{2}$ convergence preserves this norm. Since $f(0)=0$,
the limit is nonconstant. The locally uniform limit theorem for univalent
functions, a consequence of Hurwitz's theorem~\cite{DurenUnivalent},
makes $f$ univalent. Its derivative at $0$ is the limit of positive
real numbers and cannot vanish, so $f'(0)>0$.

Put $U=f(\D)$. This image is proper: otherwise $f^{-1}$ would be
a nonconstant bounded entire function, contrary to Liouville's theorem.
Its finite $\mathbf{H}^{2}$ norm gives $\E\tau_{U}<\infty$ by Lemma~\ref{pr:exit-time}.
Strong boundary convergence gives $\Re f^{*}_{j}\to\Re f^{*}$ in
$L^{2}(\T)$. Thus their distributions converge weakly. Those distributions
all equal $\mu$, so Lemma~\ref{pr:boundary-law} identifies the
real exit law of $U$ as $\mu$. Consequently $f\in\mathcal{D}_{2s,M}(\mu)$.
Applying Fatou's lemma to any coefficientwise convergent sequence
gives the displayed lower semicontinuity statement. This also proves
closedness and hence the claimed compactness. 
\end{proof}

\begin{theorem}\label{tr:minimizing} Let $0<s\leq1$ and $\mu\in\mathcal{M}_{2}$
satisfy $\mathcal{E}_{s}(U^{G}_{\mu})<\infty$. Every minimizing sequence
$U_{j}\in\mathcal{C}(\mu)$ with finite energies has a subsequence
whose normalized maps converge in $\mathcal{D}_{2s}$ to the normalized
map of a minimizing domain. If $0<s<1$, the entire sequence converges
to $H_{\mu}(z)=G_{\mu}(-z)$: 
\begin{equation}
\lim_{j\to\infty}\sum_{n\geq1}n^{2s}|c_{j,n}-(-1)^{n}a_{n}|^{2}=0.\label{tr:strong-minimizers}
\end{equation}
In particular, convergence holds in $\mathbf{H}^{2}$ and locally
uniformly in $\D$. \end{theorem} 
\begin{proof}
Set $m=\mathcal{E}_{s}(U^{G}_{\mu})$, which is the minimum by Theorem~\ref{sp:domain}.
After discarding finitely many terms, the sequence has finite, uniformly
bounded energies and therefore, by Proposition~\ref{tr:compactness},
a subsequence converging in $H^{2}$ and coefficientwise to the normalized
map $f$ of a domain $U\in\mathcal{C}(\mu)$. Lower semicontinuity
and minimality imply 
\[
m\leq\|f\|^{2}_{\mathcal{D}_{2s}}\leq\liminf_{j}\|f_{j}\|^{2}_{\mathcal{D}_{2s}}=m.
\]
The vectors $(n^{s}c_{j,n})_{n\geq1}$ are bounded in $\ell^{2}$
and converge coordinatewise to $(n^{s}c_{n})_{n\geq1}$, hence converge
weakly in $\ell^{2}$. Their squared norms converge to $m=\|f\|^{2}_{\mathcal{D}_{2s}}$.
Weak convergence together with convergence of norms in a Hilbert space
gives $f_{j}\to f$ strongly in $\mathcal{D}_{2s}$.

For $0<s<1$, the equality classification in Theorem~\ref{rg:deficit}
identifies every minimizing domain with $U^{G}_{\mu}$. Its normalized
map is $H_{\mu}$ by \eqref{pr:gross-orientation}. Apply the preceding
argument to any subsequence. Each has a further subsequence converging
in $\mathcal{D}_{2s}$ to $H_{\mu}$, which proves convergence of
the full sequence and \eqref{tr:strong-minimizers}. 
\end{proof}

\subsection{The optimal energy as a function of the law}

For $0<s\leq1$ and $\mu\in\mathcal{M}_{2}$, define the optimal value
\[
m_{s}(\mu):=\inf_{U\in\mathcal{C}(\mu)}\mathcal{E}_{s}(U)=\mathcal{E}_{s}(U^{G}_{\mu})\in(0,\infty].
\]
The second equality is Theorem~\ref{sp:domain} combined with the
main result in \cite{Energy}. Let $W_{2}$ be the quadratic Wasserstein
distance. Its one-dimensional quantile formula \cite[Proposition~2.17]{Santambrogio}
is 
\begin{equation}
W_{2}(\mu,\nu)^{2}=\inf_{\gamma\in\Pi(\mu,\nu)}\int_{\R^{2}}|x-y|^{2}\,d\gamma(x,y)=\int^{1}_{0}|q_{\mu}(t)-q_{\nu}(t)|^{2}\,dt,\label{tr:quantile-wasserstein}
\end{equation}
where $\Pi(\mu,\nu)$ denotes the couplings of $\mu$ and $\nu$.
Hardy-space stability of the Gross construction and its relation to
Wasserstein convergence are also discussed in \cite{BrownianConvergence},
cited here only for context. The following exact $p=2$ identity is
proved directly from Parseval and the quantile formula; it does not
use a result from that preprint.

\begin{proposition}\label{tr:isometry} For $\mu,\nu\in\mathcal{M}_{2}$,
with the common Gross orientation \eqref{pr:gross}, 
\begin{equation}
\|G_{\mu}-G_{\nu}\|^{2}_{\mathbf{H}^{2}}=2W_{2}(\mu,\nu)^{2}.\label{tr:exact-isometry}
\end{equation}
Consequently, for $0<r<1$, 
\begin{equation}
\sup_{|z|\leq r}|G_{\mu}(z)-G_{\nu}(z)|\leq\frac{\sqrt{2}\,r}{\sqrt{1-r^{2}}}W_{2}(\mu,\nu).\label{tr:local-bound}
\end{equation}
For every $0<s\leq1$, $m_{s}$ is lower semicontinuous under $W_{2}$
convergence to a law in $\mathcal{M}_{2}$. \end{proposition} 
\begin{proof}
The real boundary part of the difference is the even, mean-zero function
$q_{\mu}(|\theta|/\pi)-q_{\nu}(|\theta|/\pi)$. Parseval's identity
for its cosine coefficients gives 
\[
\|G_{\mu}-G_{\nu}\|^{2}_{\mathbf{H}^{2}}=2\int^{\pi}_{-\pi}|q_{\mu}(|\theta|/\pi)-q_{\nu}(|\theta|/\pi)|^{2}\,dm(\theta)=2\int^{1}_{0}|q_{\mu}-q_{\nu}|^{2}\,dt.
\]
Equation~\eqref{tr:quantile-wasserstein} proves \eqref{tr:exact-isometry}.
Cauchy--Schwarz applied to the power series, using $\sum_{n\geq1}r^{2n}=r^{2}/(1-r^{2})$,
proves \eqref{tr:local-bound}. If $\mu_{j}\to\mu$ in $W_{2}$, the
isometry gives coefficientwise convergence $a_{j,n}\to a_{n}$. Fatou's
lemma then yields 
\[
m_{s}(\mu)=\sum_{n\geq1}n^{2s}|a_{n}|^{2}\leq\liminf_{j\to\infty}\sum_{n\geq1}n^{2s}|a_{j,n}|^{2}=\liminf_{j\to\infty}m_{s}(\mu_{j}).
\]
\end{proof}

Interpolation along optimal transport plans produces Wasserstein geodesics~\cite[Theorem~5.27]{Santambrogio};
on the line this becomes linear interpolation of quantiles. We use
that representation to obtain an exact identity for the optimal spectral
energy.

\begin{theorem}\label{tr:displacement} Let $0<s\leq1$, and let
$\mu_{0},\mu_{1}\in\mathcal{M}_{2}$ have finite $m_{s}$. For $0\leq t\leq1$,
let $\mu_{t}$ have quantile 
\[
q_{t}=(1-t)q_{\mu_{0}}+tq_{\mu_{1}}.
\]
Then $\mu_{t}\in\mathcal{M}_{2}$, the curve $(\mu_{t})$ is a constant-speed
$W_{2}$ geodesic, and 
\begin{equation}
G_{\mu_{t}}=(1-t)G_{\mu_{0}}+tG_{\mu_{1}}.\label{tr:linear-gross}
\end{equation}
The optimal energy satisfies the exact identity 
\begin{align}
m_{s}(\mu_{t}) & =(1-t)m_{s}(\mu_{0})+tm_{s}(\mu_{1})-t(1-t)\|G_{\mu_{0}}-G_{\mu_{1}}\|^{2}_{\mathcal{D}_{2s}}.\label{tr:quadratic-identity}
\end{align}
In particular, 
\begin{equation}
m_{s}(\mu_{t})\leq(1-t)m_{s}(\mu_{0})+tm_{s}(\mu_{1})-2t(1-t)W_{2}(\mu_{0},\mu_{1})^{2}.\label{tr:strong-convexity}
\end{equation}
The coefficient $2$ is optimal for every $s$ in the stated range.
\end{theorem} 
\begin{proof}
The function $q_{t}$ is nondecreasing, has mean zero, and belongs
to $L^{2}(0,1)$. It is nonconstant: a convex combination with positive
coefficients of two nondecreasing functions can be constant only when
both are constant. Thus its pushforward of Lebesgue measure lies in
$\mathcal{M}_{2}$, and $q_{t}$ is its quantile almost everywhere.
The endpoint cases are immediate. The quantile formula gives 
\[
W_{2}(\mu_{t},\mu_{u})=|t-u|W_{2}(\mu_{0},\mu_{1}),\qquad0\leq t,u\leq1.
\]
Linearity of the cosine coefficients proves \eqref{tr:linear-gross}.
The Hilbert-space identity for the squared norm of a convex combination
in $\mathscr{H}_{s}$ gives \eqref{tr:quadratic-identity}. Since
$n^{2s}\geq1$, 
\[
\|G_{\mu_{0}}-G_{\mu_{1}}\|^{2}_{\mathcal{D}_{2s}}\geq\|G_{\mu_{0}}-G_{\mu_{1}}\|^{2}_{\mathbf{H}^{2}}=2W_{2}(\mu_{0},\mu_{1})^{2},
\]
which proves \eqref{tr:strong-convexity}.

For sharpness, take the centered arcsine laws on $(-R_{i},R_{i})$
with $R_{0},R_{1}>0$ and $R_{0}\ne R_{1}$. Their quantiles are $q_{\mu_{i}}(v)=-R_{i}\cos(\pi v)$
and their Gross maps are $G_{\mu_{i}}(z)=-R_{i}z$. Thus 
\[
\|G_{\mu_{0}}-G_{\mu_{1}}\|^{2}_{\mathcal{D}_{2s}}=(R_{0}-R_{1})^{2}=2W_{2}(\mu_{0},\mu_{1})^{2}.
\]
Equality holds in \eqref{tr:strong-convexity} for every $0<t<1$,
so no larger coefficient can hold uniformly. 
\end{proof}

\begin{remark}\label{tr:scope} The interpolation in Theorem~\ref{tr:displacement}
is an interpolation of quantiles and their conformal maps; no Minkowski
interpolation of the image domains is asserted. Moreover, $W_{2}$
convergence alone gives lower semicontinuity, not continuity, of positive-order
energies. Later on, remark~\ref{op:rectangle} illustrates the loss
of control of area from second moments near a degenerate limiting
law. \end{remark}

\section{Normalization and energy ratios}

\label{op:section}

For the results in this section, probability measures are centered
and nondegenerate and have a finite $p$th moment for some $p>1$,
unless a stronger moment assumption is stated. All domains on which
$\mathcal{E}_{s}$ is evaluated are simply connected and contain $0$.
We write 
\[
G_{\mu}(z)=\sum_{n\geq1}a_{n}z^{n},\qquad a_{n}=\frac{1}{\pi}\int^{\pi}_{-\pi}q_{\mu}(|\theta|/\pi)\cos(n\theta)\,d\theta,\qquad U^{G}_{\mu}=G_{\mu}(\D).
\]
In particular, we keep the increasing-quantile convention for $G_{\mu}$.
The energy $\mathcal{E}_{s}(U^{G}_{\mu})$ is allowed to be infinite.
We also use 
\[
\mathcal{E}_{0}(U)=\sum_{n\geq1}|c_{n}|^{2}
\]
 for a Riemann map $f_{U}(z)=\sum_{n\geq1}c_{n}z^{n}$ fixing $0$.
We will repeatedly use the following consequence of the Hardy-space
and exit-time identities established above: 
\begin{equation}
f_{U}\in\mathbf{H}^{2}\quad\Longrightarrow\quad\E\tau_{U}=\Var(\mu_{U})=\frac{1}{2}\mathcal{E}_{0}(U)<\infty.\label{op:h2-identities}
\end{equation}
Conversely, finiteness of $\mathcal{E}_{s}(U)$ for any $s\geq0$
implies $f_{U}\in\mathbf{H}^{2}$, because $\mathcal{E}_{0}(U)\leq\mathcal{E}_{s}(U)$.
These observations specify the integrability used below; no identity
between the exit-law variance and the mean exit time is asserted for
arbitrary weak embeddings.

\subsection{Gross normalization and the arcsine law}

We call $\mu$ \emph{Schlicht-normalized} if $R_{\mu}=1$. This convention
is equivalent to requiring $U^{G}_{\mu}$ to have conformal radius
$1$ at $0$. For a Gross-normalized measure the map $G_{\mu}(-z)$,
rather than $G_{\mu}(z)$, has derivative $+1$ at $0$ and is a normalized
Schlicht function. Both maps have the same image.

\begin{proposition} \label{op:arcsine} For every $s\geq0$, the
minimum of $\mathcal{E}_{s}(U^{G}_{\mu})$ among Schlicht -normalized
measures is $1$. The unique minimizing measure is the centered arcsine
law 
\[
d\mu_{{\rm as}}(x)=\frac{\boldsymbol{1}_{(-1,1)}(x)}{\pi\sqrt{1-x^{2}}}\,dx,
\]
and its Gross domain is $\D$. \end{proposition} 
\begin{proof}
The normalization gives $a_{1}=-1$, and therefore 
\[
\mathcal{E}_{s}(U^{G}_{\mu})=1+\sum_{n\geq2}n^{2s}|a_{n}|^{2}\geq1.
\]
Equality holds precisely when $a_{n}=0$ for every $n\geq2$. The
quantile datum is even and has mean zero, so uniqueness of Fourier
coefficients in $L^{1}(\T)$ then gives $q_{\mu}(|\theta|/\pi)=-\cos\theta$
almost everywhere. Equivalently, $q_{\mu}(t)=-\cos(\pi t)$ almost
everywhere on $(0,1)$, which is the quantile of $\mu_{{\rm as}}$.
Conversely, this law has $G_{\mu}(z)=-z$ and energy $1$ at every
exponent. 
\end{proof}

\begin{corollary} \label{op:global-schlicht} Fix $0<s\leq1$. Among
all simply connected $\mu$-domains, with $\mu$ ranging over the
Schlicht-normalized measures, the minimum of $\mathcal{E}_{s}(U)$
is $1$, and it is attained only for $U=\D$ and $\mu=\mu_{{\rm as}}$.
\end{corollary} 

\subsection{Finite-energy ratios}

In For $0<s\leq1$ and a simply connected $\mu$-domain with $\mathcal{E}_{s}(U)<\infty$,
define 
\begin{equation}
\rho_{s}(U)=\frac{\mathcal{E}_{s}(U^{G}_{\mu})}{\mathcal{E}_{s}(U)}.\label{op:ratio}
\end{equation}
The denominator is positive because $c_{1}=f_{U}'(0)\ne0$. The numerator
is positive by \eqref{pr:gross-orientation} and finite by the contraction
theorem. Thus 
\begin{equation}
0<\frac{2\Var(\mu)}{\mathcal{E}_{s}(U)}\leq\rho_{s}(U)\leq1.\label{op:ratio-bounds}
\end{equation}
Indeed, finite $\mathcal{E}_{s}(U)$ implies that $\mu$ has a finite
second moment. The Gross map then belongs to $H^{2}$, so $\mathcal{E}_{s}(U^{G}_{\mu})\geq\mathcal{E}_{0}(U^{G}_{\mu})=2\Var(\mu)$.
If $d=\operatorname{dist}(0,\partial U)$, domain monotonicity of
exit times and \eqref{op:h2-identities} also give 
\begin{equation}
\frac{d^{2}}{2}\leq\E\tau_{U}=\Var(\mu)\leq\frac{1}{2}\mathcal{E}_{s}(U^{G}_{\mu})\leq\frac{1}{2}\mathcal{E}_{s}(U).\label{op:distance-bound}
\end{equation}
Here the mean exit time from the disk of radius $d$, started at its
center, is $d^{2}/2$.

We extend \eqref{op:ratio} by the convention $\rho_{s}(U)=0$ when
$\mathcal{E}_{s}(U^{G}_{\mu})<\infty$ and $\mathcal{E}_{s}(U)=\infty$.
When both energies are infinite, $\rho_{s}(U)$ remains undefined.

The following uniqueness result shows that $U^{G}_{\mu}$ is unique
provided fulfilling some geometric constraints. The proof is in \cite{BMRemarks}.

\begin{proposition} \label{op:geometric-equality} Suppose $0<s\leq1$,
$\mathcal{E}_{s}(U)<\infty$, and $U$ is symmetric about the real
axis and vertically $\Delta$-convex: every intersection of $U$ with
a vertical line is an interval or is empty. Then $U=U^{G}_{\mu}$
and $\rho_{s}(U)=1$. \end{proposition} 

For $0<s<1$, the rigidity theorem in Section~\ref{rg:section} gives
the full characterization: under the finite-energy hypothesis, $\rho_{s}(U)=1$
if and only if $U=U^{G}_{\mu}$. Proposition~\ref{op:geometric-equality}
also supplies a sufficient geometric condition at the endpoint $s=1$.
At exponent $0$, whenever $f_{U}\in\mathbf{H}^{2}$, equation \eqref{op:h2-identities}
gives $\mathcal{E}_{0}(U)=\mathcal{E}_{0}(U^{G}_{\mu})=2\Var(\mu)$;
thus the corresponding ratio is identically $1$.

\begin{example} \label{op:uniform-example} For the uniform law $\mu$
on $(-1,1)$, the Boudabra--Markowsky domain and the Gross domain
have parametrizations 
\[
f_{{\rm BM}}(z)=-\frac{2i}{\pi}\log(1-z)=\frac{2i}{\pi}\sum_{n\geq1}\frac{z^{n}}{n},\qquad G_{\mu}(z)=-\frac{8}{\pi^{2}}\sum_{j\geq1}\frac{z^{2j-1}}{(2j-1)^{2}}.
\]
Formulas are extracted from \cite{BMNew,Gross}. Their energies satisfy
\begin{align*}
\mathcal{E}_{s}(f_{{\rm BM}}(\D)) & =\frac{4}{\pi^{2}}\sum_{n\geq1}n^{2s-2},\\
\mathcal{E}_{s}(U^{G}_{\mu}) & =\frac{64}{\pi^{4}}\sum_{j\geq1}(2j-1)^{2s-4}.
\end{align*}
For $s\geq0$, the first is finite exactly when $s<1/2$, and the
second is finite exactly when $s<3/2$. In particular, 
\[
A(f_{{\rm BM}}(\D))=\infty,\qquad A(U^{G}_{\mu})=\frac{56}{\pi^{3}}\zeta(3),
\]
and $\rho_{s}(f_{{\rm BM}}(\D))=0$ for $1/2\leq s\leq1$ under the
stated convention. \end{example} 

It is immediate that the finite-energy ratio lands in the interval
$(0,1]$. So it is natural to investigate the surjectivity of the
map 
\[
U\longmapsto\rho_{s}.
\]
For every prescribed amout of energy $c>0$, we prove surjectivity
onto $(0,1]$ when $1/2\le s\le1$. For $0<s<1/2$, a constant $\kappa_{s}>0$,
independent of the domain and the law, gives $\rho_{s}\ge\kappa_{s}$.
In other words, the area exponent $s=\frac{1}{2}$ is a threshold
for the surjectivity. This is the core result of this section. 

\subsection{Surjectivity from the area exponent onward}

\begin{theorem}[Surjectivity from the area exponent onward] \label{sur:log-surjectivity}
Let $1/2\le s\le1$ and $c>0$. Among bounded simply connected domains
$U$ containing $0$ and satisfying $\mathcal{E}_{s}(U)=c$, the range
of $\rho_{s}(U)$ is exactly $(0,1]$. The real exit law is allowed
to vary with $U$ in this statement. \end{theorem} 
\begin{proof}
First, we obtain every ratio $\eta\in(0,1]$ without prescribing the
energy to be $c$. For $0<r<1$, set 
\[
f_{r}(z)=-\frac{2i}{\pi}\log(1-rz),\qquad U_{r}=f_{r}(\D).
\]
The principal logarithm is injective on the right half-plane, so $f_{r}$
is univalent on a neighborhood of $\overline{\D}$. In particular
$U_{r}$ is bounded, contains $0$, and its real Brownian exit law
$\mu_{r}$ is the law of $\Re f_{r}(e^{i\Theta})$, where $\Theta$
is uniform on $(-\pi,\pi)$. Its coefficient energy is 
\begin{equation}
\mathcal{E}_{s}(U_{r})=\frac{4}{\pi^{2}}\sum^{\infty}_{n=1}r^{2n}n^{2s-2}.\label{sur:log-energy}
\end{equation}
We compute the quantile of $\mu_{r}$. Put $A_{r}=\arg(1-re^{i\Theta})$.
This random variable is symmetric and is supported on $[-\arcsin r,\arcsin r]$.
If $0<x<\arcsin r$, then 
\[
\{A_{r}>x\}=\{r\sin(x-\Theta)>\sin x\},
\]
because $\Re(1-re^{i\Theta})>0$. Consequently 
\[
\mathbb{P}(A_{r}>x)=\frac{1}{2}-\frac{1}{\pi}\arcsin\frac{\sin x}{r}.
\]
Using symmetry and $\Re f_{r}(e^{i\Theta})=2A_{r}/\pi$, we obtain
\begin{equation}
q_{r}(t)=-\frac{2}{\pi}\arcsin(r\cos(\pi t)),\qquad0<t<1.\label{sur:log-quantile}
\end{equation}
In particular 
\[
0\le q_{r}'(t)=\frac{2r\sin(\pi t)}{\sqrt{1-r^{2}\cos^{2}(\pi t)}}\le2.
\]
For the Gross boundary datum $q_{r}(|\theta|/\pi)$, Parseval's identity
therefore gives 
\begin{equation}
\mathcal{E}_{s}(U^{G}_{\mu_{r}})\le\mathcal{E}_{1}(U^{G}_{\mu_{r}})=\frac{2}{\pi^{2}}\int^{1}_{0}|q_{r}'(t)|^{2}\,dt\le\frac{8}{\pi^{2}}.\label{sur:log-numerator-bound}
\end{equation}
Since the series in \eqref{sur:log-energy} diverges as $r\uparrow1$
whenever $s\ge1/2$, equations \eqref{sur:log-energy} and \eqref{sur:log-numerator-bound}
imply 
\[
\rho_{s}(U_{r})\longrightarrow0\qquad(r\uparrow1).
\]
For every compact subinterval of $(0,1)$, formula \eqref{sur:log-quantile}
varies continuously in $H^{1}(0,1)$. The associated Gross coefficients
therefore vary continuously in the $\mathcal{E}_{1}$ norm, hence
also in the $\mathcal{E}_{s}$ norm. The denominator in \eqref{sur:log-energy}
is positive and continuous. Thus $r\mapsto\rho_{s}(U_{r})$ is continuous
on $(0,1)$. As $r\downarrow0$, 
\[
\frac{q_{r}(t)}{r}\longrightarrow-\frac{2}{\pi}\cos(\pi t)\quad\hbox{in}H^{1}(0,1),\qquad\frac{\mathcal{E}_{s}(U_{r})}{r^{2}}\longrightarrow\frac{4}{\pi^{2}}.
\]
The first limit gives $\mathcal{E}_{s}(U^{G}_{\mu_{r}})/r^{2}\to4/\pi^{2}$,
so $\rho_{s}(U_{r})\to1$. The intermediate value theorem yields every
value in $(0,1)$, and a centered disk yields $1$.

Now, we finalize the proof with prescribed energy. For a domain $U$
with finite positive energy, let $\lambda=(c/\mathcal{E}_{s}(U))^{1/2}$.
Positive dilation multiplies both energies by $\lambda^{2}$: its
real exit law has quantile $\lambda q$, whose Gross map is $\lambda G_{\mu}$.
Hence $\lambda U$ has energy $c$ and the same ratio as $U$. The
universal bounds $0<\rho_{s}(U)\le1$ complete the proof. 
\end{proof}

\begin{remark} The same construction does not prove surjectivity
for $0<s<1/2$. In that range its denominator tends to $4\zeta(2-2s)/\pi^{2}<\infty$,
whereas its Gross numerator tends to the strictly positive energy
of the uniform law on $(-1,1)$. \end{remark}

\subsection{A positive lower bound below the area exponent}

\label{sur:lower-section}

The behavior of the ratio changes at $s=1/2$. We prove that the surjectivity
statement in Theorem~\ref{sur:log-surjectivity} cannot extend to
$0<s<1/2$. The constants below depend only on $s$; no assertion
of sharpness is made. As usual, the circle has normalized measure
$dm(\theta)=d\theta/(2\pi)$. The constant $c_{s}$ is again $\frac{\Gamma(1+2s)\sin(\pi s)}{\pi}$.
For a real measurable function $u$ on the circle, let $h$ denote
the decreasing rearrangement of $|u|$ on $(0,1)$, and set 
\begin{equation}
L_{s}(u)=\int^{1}_{0}t^{-2s}h(t)^{2}\,dt.\label{sur:lorentz-functional}
\end{equation}
Thus $L_{s}(u)^{1/2}$ is the usual rearrangement expression for the
Lorentz space $L^{p,2}$, where $p=2/(1-2s)$; see \cite[Example~1.1, p.~284]{Peetre}.
Only the explicit integral in \eqref{sur:lorentz-functional} is used
in the proof. The Lorentz refinement of fractional Sobolev embedding
is classical; see Peetre~\cite[Section~8]{Peetre}. The following
proof gives the periodic, mean-zero estimate and an explicit constant
in our normalization.

\begin{lemma}[A fractional Hardy estimate and a Lorentz bound]
\label{sur:hardy} Let $0<s<1/2$, put $p=2/(1-2s)$, and use normalized
measure $d\theta/(2\pi)$ on $\T$. For a real mean-zero function
$u\in H^{s}(\T)$, let $h:(0,1)\to[0,\infty)$ be the decreasing rearrangement
of $|u|$. Then 
\begin{equation}
\|u\|^{2}_{L^{p,2}}:=\int^{1}_{0}h(t)^{2}t^{-2s}\,dt\le B_{s}\mathscr{S}_{s}(u),\label{sur:lorentz}
\end{equation}
where $B_{s}<\infty$ depends only on $s$. For example, with the
normalization $\mathscr{S}_{s}(u)=c_{s}Q_{s}(u)/(4\pi)$ used above,
one may take 
\[
B_{s}=\frac{2^{4+2s}\pi^{2s}/c_{s}+2\,3^{2s}}{(1-2^{s-1/2})^{2}}.
\]
No sharpness of this constant is asserted. \end{lemma} 
\begin{proof}
We first establish a weighted estimate for a bounded nonnegative nonincreasing
function $w$ on $(0,\pi)$. Extend $w$ evenly to the circle, and
put 
\[
L=\frac{\pi}{3},\qquad H^{2}=\int^{\pi}_{0}w(x)^{2}x^{-2s}\,dx,\qquad D=\int^{L}_{0}\bigl(w(x)-w(2x)\bigr)^{2}x^{-2s}\,dx.
\]
Since $2s<1$ and $w$ is bounded, $H<\infty$. In $L^{2}((0,\pi),x^{-2s}dx)$,
decompose $w$ as the sum of $(w(x)-w(2x))\mathbf{1}_{(0,L)}(x)$,
$w(2x)\mathbf{1}_{(0,L)}(x)$, and $w(x)\mathbf{1}_{[L,\pi)}(x)$.
Minkowski's inequality and a change of variables give 
\[
H\le\sqrt{D}+2^{s-1/2}H+L^{-s}\|w\|_{L^{2}(0,\pi)}.
\]
For $0<x<L$ and $2x<y<3x$, monotonicity gives $w(x)-w(y)\ge w(x)-w(2x)\ge0$.
Consequently, 
\begin{align*}
Q_{s}(w) & \ge\int^{L}_{0}\int^{3x}_{2x}\frac{|w(x)-w(y)|^{2}}{|x-y|^{1+2s}}\,dy\,dx\\
 & \ge2^{-1-2s}D.
\end{align*}
Here we used $K_{s}(x-y)\ge|x-y|^{-1-2s}$, and restricted the nonnegative
double integral defining $Q_{s}$ to the indicated region. It follows
that 
\begin{equation}
\int^{\pi}_{0}w(x)^{2}x^{-2s}\,dx\le\frac{2}{(1-2^{s-1/2})^{2}}\left(2^{1+2s}Q_{s}(w)+\left(\frac{3}{\pi}\right)^{2s}\|w\|^{2}_{L^{2}(0,\pi)}\right).\label{sur:weighted-hardy}
\end{equation}
Now set $w(x)=h(|x|/\pi)$ on $\T$. Circular rearrangement and $\bigl||a|-|b|\bigr|\le|a-b|$
imply 
\[
Q_{s}(w)\le Q_{s}(|u|)\le Q_{s}(u),\qquad\|w\|^{2}_{L^{2}(0,\pi)}=\pi\|u\|^{2}_{L^{2}(\T)}.
\]
If $w$ is unbounded, apply \eqref{sur:weighted-hardy} to $w_{M}=\min(w,M)$.
These truncations are still nonnegative and nonincreasing, and satisfy
$Q_{s}(w_{M})\le Q_{s}(w)$ and $\|w_{M}\|_{2}\le\|w\|_{2}$. Monotone
convergence as $M\to\infty$ therefore proves \eqref{sur:weighted-hardy}
for $w$ as well. Finally, 
\[
\int^{1}_{0}h(t)^{2}t^{-2s}\,dt=\pi^{2s-1}\int^{\pi}_{0}w(x)^{2}x^{-2s}\,dx.
\]
Insert \eqref{sur:weighted-hardy}, use $Q_{s}(u)=4\pi\mathscr{S}_{s}(u)/c_{s}$,
and apply the mean-zero Fourier Poincare inequality $\|u\|^{2}_{L^{2}(\T)}\le\mathscr{S}_{s}(u)$.
This gives \eqref{sur:lorentz} with the displayed constant. 
\end{proof}

Maximum-modulus characterizations of univalent functions in weighted
Dirichlet spaces are established by Qian and Shi~\cite[Main Theorem, pp.~1214--1217]{QianShi}.
The choice $K(t)=t^{1-2s}$ gives the power-weight framework relevant
here. We include the weighted-area argument and then estimate the
maximum modulus from the real boundary datum.

\begin{lemma}[Univalent reverse estimate]\label{sur:univalent-reverse}
Suppose $0<s<1/2$, $f\in\mathbf{H}^{2}$ is univalent on $\D$, and
$f(0)=0$. If $u=\Re f^{*}$, then 
\[
\|f\|^{2}_{\mathcal{D}_{2s}}=\sum_{n\ge1}n^{2s}|c_{n}|^{2}\le A_{s}L_{s}(u),\qquad f(z)=\sum_{n\ge1}c_{n}z^{n},
\]
for a finite constant $A_{s}$ independent of $f$. The assertion
is understood in the extended nonnegative real numbers. \end{lemma} 
\begin{proof}
In this proof, $C$ denotes an absolute positive constant and $C_{s}$
a positive constant depending only on $s$; their values may change
between occurrences. Put $\alpha=1-2s\in(0,1)$ and $M_{f}(r)=\max_{|z|=r}|f(z)|$.
Parseval and Tonelli give the exact identity 
\begin{equation}
\sum_{n\ge1}n^{2s}|c_{n}|^{2}=\frac{2^{\alpha}}{\pi\Gamma(1+\alpha)}\int_{\D}|f'(z)|^{2}\bigl(\log(1/|z|)\bigr)^{\alpha}\,dA(z),\label{sur:weighted-area-identity}
\end{equation}
where $dA$ is ordinary planar area. Indeed, the contribution of the
$n$th coefficient to the integral is 
\[
2\pi n^{2}|c_{n}|^{2}\int^{1}_{0}r^{2n-1}(\log(1/r))^{\alpha}\,dr=\pi2^{-\alpha}\Gamma(1+\alpha)n^{1-\alpha}|c_{n}|^{2}.
\]
Write $B_{f}(r)=\int_{|z|<r}|f'(z)|^{2}\,dA(z)$. Univalence and the
maximum principle imply 
\begin{equation}
B_{f}(r)=\Area(f(r\D))\le\pi M_{f}(r)^{2}.\label{sur:univalent-area-bound}
\end{equation}
This is the step that uses injectivity: otherwise the integral counts
area with multiplicity. Using 
\[
(\log(1/\rho))^{\alpha}=\alpha\int^{1}_{\rho}\frac{(\log(1/r))^{\alpha-1}}{r}\,dr
\]
and Tonelli, the integral in \eqref{sur:weighted-area-identity} equals
\[
\alpha\int^{1}_{0}B_{f}(r)\frac{(\log(1/r))^{\alpha-1}}{r}\,dr.
\]
For $1/2\le r<1$, its weight is bounded by a constant depending only
on $\alpha$ times $(1-r)^{\alpha-1}$. For $r\le1/2$, Schwarz's
lemma applied to $f(z)/M_{f}(1/2)$ on the disk of radius $1/2$ gives
$M_{f}(r)\le2rM_{f}(1/2)$. Thus the contribution of $r<1/2$ is at
most a constant times $M_{f}(1/2)^{2}$. This last quantity is controlled
by the integral over $1/2\le r\le3/4$, since $M_{f}$ is nondecreasing.
Together with \eqref{sur:univalent-area-bound}, this proves 
\begin{equation}
\sum_{n\ge1}n^{2s}|c_{n}|^{2}\le C_{s}\int^{1}_{0}h^{\alpha-1}M_{f}(1-h)^{2}\,dh.\label{sur:maximum-modulus-bound}
\end{equation}

We bound the last integral directly from the real boundary values,
using the Schwarz integral formula~\cite{DurenHardy}. This formula,
valid because $f\in H^{2}$ and $f(0)=0$, is 
\[
f(z)=\int^{\pi}_{-\pi}\frac{e^{i\theta}+z}{e^{i\theta}-z}\,u(\theta)\,dm(\theta).
\]
For $|z|=1-h\ge1/2$, the absolute value of its kernel is bounded
by $C/(h+\operatorname{dist}_{\T}(\theta,\arg z))$. Its decreasing
rearrangement is consequently at most $C/(h+t)$ on $(0,1)$. The
Hardy--Littlewood rearrangement inequality \cite[Chapter~2, Section~2]{BennettSharpley}
therefore gives, with $V=u^{\mathrm{abs}}$, 
\begin{equation}
M_{f}(1-h)\le C\left(\frac{1}{h}\int^{h}_{0}V(t)\,dt+\int^{1}_{h}\frac{V(t)}{t}\,dt\right).\label{sur:schwarz-rearrangement}
\end{equation}
The same bound holds for $h\ge1/2$: the Schwarz kernel is then uniformly
bounded, and the first term on the right dominates $\int^{1}_{0}V(t)\,dt$
because $V$ is nonincreasing.

For completeness, both weighted Hardy estimates required for \eqref{sur:schwarz-rearrangement}
follow from Minkowski's inequality. Extend $V$ by zero to $(1,\infty)$
and use the measure $h^{\alpha-1}dh$ on $(0,\infty)$. Since 
\[
\frac{1}{h}\int^{h}_{0}V(t)\,dt=\int^{1}_{0}V(hv)\,dv,\qquad\int^{1}_{h}\frac{V(t)}{t}\,dt=\int^{\infty}_{1}V(hv)\,\frac{dv}{v}\quad(0<h<1),
\]
the corresponding operator norms into $L^{2}(h^{\alpha-1}dh)$ are
at most 
\[
\int^{1}_{0}v^{-\alpha/2}\,dv=\frac{2}{2-\alpha},\qquad\int^{\infty}_{1}v^{-1-\alpha/2}\,dv=\frac{2}{\alpha}.
\]
Consequently 
\[
\int^{1}_{0}h^{\alpha-1}M_{f}(1-h)^{2}\,dh\le C_{s}\int^{1}_{0}t^{\alpha-1}V(t)^{2}\,dt=C_{s}L_{s}(u).
\]
Combining this with \eqref{sur:maximum-modulus-bound} proves the
lemma. If the right-hand side is infinite the conclusion is immediate;
otherwise the preceding estimates also establish finiteness of the
left-hand side. 
\end{proof}

We combine this reverse estimate with the fractional rearrangement
bound from Lemma~\ref{sur:hardy}: for every real mean-zero circle
function $v$ of finite $\mathscr{S}_{s}$ energy, 
\begin{equation}
L_{s}(v)\le B_{s}\mathscr{S}_{s}(v),\qquad0<s<1/2.\label{sur:gross-lorentz}
\end{equation}
The periodic estimate, with its normalization, was proved in Lemma~\ref{sur:hardy}.

\begin{remark} contributes something essential. For the nonunivalent
functions 
\[
f_{N}(z)=z^{N},\,\,\,\,\,\,\,\,u_{N}=\cos(N\theta).
\]
The distributions of $|u_{N}|$ are identical and so $L_{s}(u_{N}0$
is independent of $N$. Nevertheless, 
\[
\sum_{n\ge1}n^{2s}|c_{n}|^{2}=N^{2s}\longrightarrow\infty.
\]

\end{remark}

\begin{theorem}[The ratio stays positive below the area exponent]
\label{sur:positive-infimum} For every $0<s<1/2$, there exists $\kappa_{s}>0$
depending only on $s$ such that every simply connected $U\ni0$ with
$\mathcal{E}_{s}(U)<\infty$ satisfies 
\begin{equation}
\kappa_{s}\mathcal{E}_{s}(U)\le\mathcal{E}_{s}(U^{G}_{\mu_{U}})\le\mathcal{E}_{s}(U).\label{sur:two-sided-energy}
\end{equation}
In particular, at any prescribed positive value of $\mathcal{E}_{s}(U)$,
the ratio $\rho_{s}(U)$ is bounded below by the same $\kappa_{s}$,
and its range is not $(0,1]$. \end{theorem} 
\begin{proof}
Let $f:\D\to U$ be a conformal map with $f(0)=0$, and put $u=\Re f^{*}$.
Finite energy implies $f\in\mathbf{H}^{2}$ and identifies the law
of $u$ with $\mu_{U}$. If $q$ is its quantile and $v(\theta)=q(|\theta|/\pi)$,
then $u$ and $v$ are equimeasurable. Consequently $L_{s}(u)=L_{s}(v)$.
Theorem~\ref{sp:domain} gives finite $\mathcal{E}_{s}(U^{G}_{\mu_{U}})$,
so \eqref{sur:gross-lorentz} applies to $v$. Lemma~\ref{sur:univalent-reverse}
and \eqref{sur:gross-lorentz} give 
\[
\mathcal{E}_{s}(U)\le A_{s}L_{s}(u)=A_{s}L_{s}(v)\le\frac{A_{s}B_{s}}{2}\mathcal{E}_{s}(U^{G}_{\mu_{U}}).
\]
Taking $\kappa_{s}=\min\{1,2/(A_{s}B_{s})\}$ proves the lower bound.
The upper bound is the spectral rearrangement theorem already proved.
Dividing by the positive denominator gives the assertion about $\rho_{s}$. 
\end{proof}

\begin{remark} This theorem establishes a strictly positive universal
lower endpoint for $0<s<1/2$ but does not determine its sharp value
or whether that value is attained. Combined with Theorem~\ref{sur:log-surjectivity},
it identifies $s=1/2$ as the threshold for finite-energy ratios tending
to zero. At $s=0$, the finite-$H^{2}$ ratio is identically $1$.

The finite-energy restriction is essential. For the weak domains $V_{a,b}$
of Proposition~\ref{pr:weak-counterexample}, \eqref{pr:exit-coefficients}
gives $\mathcal{E}_{0}(V_{a,b})=\infty$, and hence $\mathcal{E}_{s}(V_{a,b})=\infty$
for every $s\geq0$. Their balanced Bernoulli exit law has finite
Gross energy for $0<s<1/2$ by Theorem~\ref{cr:diameter}. Thus the
extended convention of Section~\ref{op:section} gives $\rho_{s}(V_{a,b})=0$
in that range; these domains do not satisfy the finite-energy hypothesis
above. \end{remark}

The restriction$s<\tfrac{1}{2}$ in the lower-bound argument is reflected
in the Sobolev--Lorentz relation
\[
\frac{1}{p}=\frac{1}{2}-s,\qquad L_{s}(u)=\int^{1}_{0}\bigl[t^{1/p}u^{\mathrm{abs}}(t)\bigr]^{2}\,\frac{dt}{t}.
\]
For $0<s<\tfrac{1}{2}$, the exponent $p$ is finite, and the preceding
Lorentz estimate controls $L_{s}(u)$ by $\mathscr{S}_{s}(u)$ for
mean-zero $u$. As $s\uparrow\tfrac{1}{2}$, however, $p\to\infty$,
and the formal endpoint functional becomes
\[
L_{1/2}(u)=\int^{1}_{0}\frac{|u^{\mathrm{abs}}(t)|^{2}}{t}\,dt.
\]
This integral is infinite for every nonzero $u$. Indeed, there exist
$a,\delta>0$ such that $u^{\mathrm{abs}}(t)\ge a$ for $0<t<\delta$,
and therefore
\[
L_{1/2}(u)\ge a^{2}\int^{\delta}_{0}\frac{dt}{t}=\infty.
\]
The same restriction appears in the univalent reverse estimate through
the parameter $\alpha=1-2s$: the weighted Hardy bound $2/\alpha$
requires $\alpha>0$. These observations explain the range of the
particular estimates used here, but their degeneration alone would
not establish a sharp threshold for the energy ratio. The threshold
itself is verified by the logarithmic family
\[
f_{r}(z)=-\frac{2i}{\pi}\log(1-rz),\qquad U_{r}=f_{r}(\mathbb{D}),\qquad0<r<1.
\]
Each $U_{r}$ is bounded and has finite energy, while
\[
\mathcal{E}_{s}(U_{r})=\frac{4}{\pi^{2}}\sum_{n\ge1}r^{2n}n^{2s-2}\longrightarrow\infty\quad\text{as }r\uparrow1
\]
precisely when $s\ge\tfrac{1}{2}$. In contrast, the corresponding
canonical energies satisfy
\[
\mathcal{E}_{s}(U^{G}_{\mu_{r}})\le\frac{8}{\pi^{2}},\qquad0<s\le1.
\]
Consequently, $\rho_{s}(U_{r})\to0$ for $\tfrac{1}{2}\le s\le1$.
Together with the positive universal lower bound for $0<s<\tfrac{1}{2}$,
this identifies $\tfrac{1}{2}$ as the exact threshold for finite-energy
ratios approaching zero within $0<s\le1$. Positive dilation preserves
the ratio, so the same threshold holds at every prescribed positive
energy. Its geometric significance is given by
\[
\operatorname{Area}(U)=\pi\mathcal{E}_{1/2}(U):
\]
the logarithmic family develops unbounded area while its canonical
counterparts retain uniformly bounded area.

\section{Sharp bounds and a criterion for finite area}

\label{cr:section}

The coefficient energies can also be estimated directly from the exit
law. Throughout this section, $\mu$ is centered and nondegenerate,
and its essential support endpoints 
\[
a=\mathop{{\rm ess\,inf}}_{\xi\sim\mu}x,\qquad b=\mathop{{\rm ess\,sup}}_{\xi\sim\mu}x
\]
are finite. Write $D=b-a>0$, $q=q_{\mu}$, and let $\nu=dq$ be the
Stieltjes measure, or equivalently the distributional derivative,
of the nondecreasing quantile on $(0,1)$. Thus $\nu$ is a finite
positive measure of total mass $D$. We use the endpoint limits $q(0^{+})=a$
and $q(1^{-})=b$; no endpoint atoms are added to $\nu$. Stieltjes
integration by parts gives 
\begin{equation}
a_{n}=2\int^{1}_{0}q(t)\cos(\pi nt)\,dt=-\frac{2}{\pi n}\int^{1}_{0}\sin(\pi nt)\,d\nu(t).\label{cr:stieltjes}
\end{equation}
This identity is valid whether or not $q$ is absolutely continuous.

\begin{theorem}[Sharp bound in terms of support diameter] \label{cr:diameter}
For every $0<s<1/2$, 
\begin{equation}
\inf_{U\in\mathcal{C}(\mu)}\mathcal{E}_{s}(U)=\mathcal{E}_{s}(U^{G}_{\mu})\leq\frac{4D^{2}}{\pi^{2}}\bigl(1-2^{-(2-2s)}\bigr)\zeta(2-2s).\label{cr:diameter-bound}
\end{equation}
Equality in the upper bound holds if and only if $\mu=\tfrac{1}{2}(\delta_{-D/2}+\delta_{D/2})$.
In particular, every compactly supported law in the stated class admits
an embedding with finite $\mathcal{E}_{s}$ for $s<1/2$. \end{theorem} 
\begin{proof}
The minimum identity follows from Theorem~\ref{sp:domain}. Put $\alpha=2-2s>1$
and define 
\[
S_{\alpha}(t)=\sum_{n\geq1}\frac{\sin^{2}(\pi nt)}{n^{\alpha}},\qquad0\leq t\leq1.
\]
The series converges uniformly. We first show that its unique maximum
is at $t=1/2$. Set $C_{\alpha}(x)=\sum_{n\geq1}n^{-\alpha}\cos(nx)$,
so that $S_{\alpha}(t)=\tfrac{1}{2}(\zeta(\alpha)-C_{\alpha}(2\pi t))$.
The gamma integral and absolute summability give 
\[
C_{\alpha}(x)=\frac{1}{\Gamma(\alpha)}\int^{\infty}_{0}u^{\alpha-1}\frac{r\cos x-r^{2}}{1-2r\cos x+r^{2}}\,du,\qquad r=e^{-u}.
\]
For fixed $0<r<1$, the derivative of the fraction with respect to
$c=\cos x$ is $r(1-r^{2})/(1-2rc+r^{2})^{2}>0$. Consequently $C_{\alpha}$
is strictly decreasing on $[0,\pi]$. Its symmetry about $\pi$ proves
the asserted unique maximum of $S_{\alpha}$, with 
\[
S_{\alpha}(1/2)=\sum_{j\geq1}(2j-1)^{-\alpha}=(1-2^{-\alpha})\zeta(\alpha).
\]
Apply Cauchy--Schwarz with respect to $\nu$ in \eqref{cr:stieltjes},
and then Tonelli's theorem: 
\begin{align*}
\mathcal{E}_{s}(U^{G}_{\mu}) & =\frac{4}{\pi^{2}}\sum_{n\geq1}\frac{1}{n^{\alpha}}\left(\int\sin(\pi nt)\,d\nu(t)\right)^{2}\\
 & \leq\frac{4D}{\pi^{2}}\int S_{\alpha}(t)\,d\nu(t)\leq\frac{4D^{2}}{\pi^{2}}S_{\alpha}(1/2).
\end{align*}
If equality holds in the final bound, the unique maximum just proved
forces $\nu=D\delta_{1/2}$. Thus $q$ equals a constant $a$ on $(0,1/2)$
and the constant $a+D$ on $(1/2,1)$. Centering gives $a=-D/2$,
which identifies $\mu$. Conversely, this two-point law has $\nu=D\delta_{1/2}$
and makes both inequalities equalities. 
\end{proof}

\begin{remark}\label{cr:endpoint} At $s=0$, the same calculation
yields $2\operatorname{Var}(\mu)=\mathcal{E}_{0}(U^{G}_{\mu})\leq D^{2}/2$,
which is an independent Fourier--Stieltjes derivation of Popoviciu’s
sharp variance bound for an interval of length $D$; see~\cite{Popoviciu}.
The threshold $s<1/2$ in Theorem~\ref{cr:diameter} cannot be increased:
for its two-point extremizer, 
\[
\mathcal{E}_{s}(U^{G}_{\mu})=\frac{4D^{2}}{\pi^{2}}\sum_{j\geq1}(2j-1)^{2s-2},
\]
which is infinite for $s\geq1/2$. \end{remark}

\begin{theorem}[A law-based criterion for finite area] \label{cr:area}
Define the nonnegative kernel on $(0,1)^{2}$ by 
\[
L(t,u)=\log\frac{\sin(\pi(t+u)/2)}{\sin(\pi|t-u|/2)},\qquad L(t,t)=+\infty.
\]
Then, with extended nonnegative values, 
\begin{equation}
\inf_{U\in\mathcal{C}(\mu)}\operatorname{Area}(U)=\operatorname{Area}(U^{G}_{\mu})=\frac{2}{\pi}\iint_{(0,1)^{2}}L(t,u)\,d\nu(t)\,d\nu(u).\label{cr:area-formula}
\end{equation}
Thus a simply connected embedding of finite area exists if and only
if the integral in \eqref{cr:area-formula} is finite. \end{theorem} 
\begin{proof}
For $0<r<1$, the absolutely and uniformly convergent series 
\[
B_{r}(t,u)=\sum_{n\geq1}\frac{r^{n}}{n}\sin(\pi nt)\sin(\pi nu)
\]
can be integrated term by term against the finite measure $\nu\otimes\nu$.
Summing the logarithmic power series gives 
\begin{equation}
B_{r}(t,u)=\frac{1}{4}\log\frac{1-2r\cos(\pi(t+u))+r^{2}}{1-2r\cos(\pi(t-u))+r^{2}}.\label{cr:abel}
\end{equation}
The formula also proves positivity and monotonicity in $r$, which
are not immediate from the sine series. Indeed, write 
\[
A=4\sin^{2}(\pi(t-u)/2),\quad B=4\sin^{2}(\pi(t+u)/2),\quad h(r)=r+r^{-1}-2.
\]
Here $B-A=4\sin(\pi t)\sin(\pi u)>0$, and the ratio in \eqref{cr:abel}
is $(h(r)+B)/(h(r)+A)$. Since $h$ decreases on $(0,1)$, this ratio
increases and is at least one. Hence $B_{r}(t,u)\uparrow\tfrac{1}{2}L(t,u)$
as $r\uparrow1$, including the infinite diagonal limit.

Using \eqref{cr:stieltjes} at each $r<1$ yields 
\[
\sum_{n\geq1}nr^{n}a^{2}_{n}=\frac{4}{\pi^{2}}\iint B_{r}(t,u)\,d\nu(t)\,d\nu(u).
\]
Monotone convergence on both sides now gives $\mathcal{E}_{1/2}(U^{G}_{\mu})=\frac{2}{\pi^{2}}\iint L\,d\nu\,d\nu$.
Multiply by $\pi$ and apply Theorem~\ref{sp:domain} to obtain \eqref{cr:area-formula}.
If the integral is finite, the Gross domain supplies a finite-area
embedding; if it is infinite, the minimum identity excludes every
finite-area competitor. 
\end{proof}

\begin{corollary}[Quantile concentration and support gaps] \label{cr:gaps}
Let $I\subset[\delta,1-\delta]$ be an interval of length $\varepsilon>0$,
where $0<\delta<1/2$, and put $m=\nu(I)$. Every $U\in\mathcal{C}(\mu)$
satisfies 
\begin{equation}
\operatorname{Area}(U)\geq\frac{2m^{2}}{\pi}\left[\log\frac{2\sin(\pi\delta)}{\pi\varepsilon}\right]_{+},\qquad[x]_{+}=\max\{x,0\}.\label{cr:concentration}
\end{equation}
If the support of $\mu$ has an internal gap, every simply connected
$\mu$-domain has infinite area. In particular, this holds for every
nondegenerate finitely supported law. \end{corollary} 
\begin{proof}
For $t,u\in I$, the numerator in $L(t,u)$ is at least $\sin(\pi\delta)$,
whereas its denominator is at most $\pi\varepsilon/2$. Since $L\geq0$,
restricting \eqref{cr:area-formula} to $I\times I$ proves \eqref{cr:concentration}.

Suppose $(x_{-},x_{+})$ is a support gap with positive mass on both
sides. Then $q$ has a jump of size at least $x_{+}-x_{-}>0$ at the
rank $p=\mu(({-\infty},x_{-}])\in(0,1)$. Equivalently, $\nu(\{p\})>0$.
The diagonal value $L(p,p)=+\infty$ implies that the integral in
\eqref{cr:area-formula} is infinite. 
\end{proof}

The support-gap obstruction is consistent with the geometric observation
in \cite{Gross}: a gap forces a vertical strip in the Gross domain.
Formula~\eqref{cr:concentration} also quantifies the area cost when
a fixed quantile increase occurs across a small, positive interval
of probability ranks. An atom of $\mu$ instead produces a constant
interval of $q$; it should not be confused with an atom of the derivative
measure $\nu$, which records a support gap. In \cite{BoundedGross},
the authors provide a full characterization of the boundedness of
Gross domains when $\mu$ has a density. For the sake of completeness,
we illustrate an example of a distribution $\mu$ that generates an
unbounded Gross domain where we apply our approach.
\begin{example}[Connected support does not ensure finite area]
\label{cr:smooth-obstruction} Set $B=(\log2)^{-1/2}$ and define
\[
\rho(x)=\begin{cases}
2|x|^{-3}e^{-1/x^{2}}, & 0<|x|<B,\\
0, & x=0\text{ or }|x|\geq B.
\end{cases}
\]
Then $d\mu(x)=\rho(x)\,dx$ is a centered probability measure with
connected support $[-B,B]$. Its density is smooth on the interior
$(-B,B)$ and positive except at $0$, yet every $U\in\mathcal{C}(\mu)$
has infinite area. 
\end{example}

\begin{proof}
For $x>0$, the derivative of $e^{-1/x^{2}}$ is $2x^{-3}e^{-1/x^{2}}$.
Hence the integral of $\rho$ over $(0,B)$ is $e^{-1/B^{2}}=1/2$.
Symmetry gives total mass one and mean zero. All derivatives of $\rho$
vanish at $0$, so the density is smooth through that point; no smoothness
at $\pm B$ is asserted. Since $\rho>0$ elsewhere in $(-B,B)$, its
support is the full interval $[-B,B]$.For $0<|x|<B$ the distribution
function is 
\[
F_{\mu}(x)=\frac{1}{2}+\operatorname{sgn}(x)e^{-1/x^{2}}.
\]
Writing $p=1/2$, inversion gives the continuous quantile 
\[
q(t)=\frac{\operatorname{sgn}(t-p)}{\sqrt{\log(1/|t-p|)}}\quad(t\ne p),\qquad q(p)=0.
\]
Its Stieltjes measure therefore satisfies 
\begin{equation}
\nu((p-r,p+r))=\frac{2}{\sqrt{\log(1/r)}},\qquad0<r<1/2.\label{cr:smooth-mass}
\end{equation}
Let $J=[1/4,3/4]$ and $c=\sqrt{2}/\pi$. For $t,u\in J$, the numerator
of $L(t,u)$ is at least $\sin(\pi/4)$ and its denominator is at
most $\pi|t-u|/2$. Together with $L\geq0$, this implies 
\[
L(t,u)\geq\left[\log\frac{c}{|t-u|}\right]_{+}.
\]
The nonnegative layer-cake formula and Tonelli's theorem yield 
\begin{align*}
\iint L(t,u)\,d\nu(t)\,d\nu(u) & \geq\int^{\infty}_{0}(\nu\otimes\nu)\{(t,u)\in J^{2}:|t-u|<ce^{-v}\}\,dv\\
 & \geq\int^{\infty}_{0}\nu(I_{v})^{2}\,dv,\qquad I_{v}=(p-\tfrac{c}{2}e^{-v},p+\tfrac{c}{2}e^{-v}).
\end{align*}
Indeed, $c/2<1/4$, so $I_{v}\subset J$, and every pair of its points
has distance less than $ce^{-v}$. By \eqref{cr:smooth-mass}, the
last integral equals 
\[
4\int^{\infty}_{0}\frac{dv}{v+\log(2/c)}=+\infty.
\]
Theorem~\ref{cr:area} now excludes every simply connected embedding
of finite area for this law. 
\end{proof}

\section{Fixed area and the exit-law variance}

\label{op:fixed-area-section}

For $c>0$, let 
\[
\mathcal{A}_{c}=\{U\subset\C:\ U\text{ is simply connected},\ 0\in U,\ A(U)=c\}.
\]
Every $U\in\mathcal{A}_{c}$ has $f_{U}\in\mathbf{H}^{2}$, since
$\mathcal{E}_{0}(U)\leq\mathcal{E}_{1/2}(U)=c/\pi$. Its real exit
law is therefore centered and has finite, positive variance.

\begin{theorem} \label{op:variance-range} For every $c>0$, 
\[
\{\Var(\mu_{U}):U\in\mathcal{A}_{c}\}=\left(0,\frac{c}{2\pi}\right].
\]
The upper endpoint is attained only by the disk centered at $0$ with
area $c$. The infimum is $0$ and is not attained. \end{theorem} 
\begin{proof}
Equation~\eqref{op:h2-identities} and the area formula yield 
\[
0<\Var(\mu_{U})=\frac{1}{2}\sum_{n\geq1}|c_{n}|^{2}\leq\frac{1}{2}\sum_{n\geq1}n|c_{n}|^{2}=\frac{c}{2\pi}.
\]
Strict positivity follows from $c_{1}\ne0$. Equality in the upper
bound holds if and only if $c_{n}=0$ for $n\geq2$, which characterizes
the stated disk.

To obtain every intermediate value, put $R=\sqrt{c/\pi}$ and, for
$0\leq\kappa<1$, let 
\[
U_{\kappa}=\{w\in\C:|w-i\kappa R|<R\}.
\]
These disks contain $0$ and have area $c$. Their Riemann maps fixing
$0$ and having positive derivative there are 
\begin{align*}
f_{\kappa}(z) & =R\left(\frac{z-i\kappa}{1+i\kappa z}+i\kappa\right)=R\frac{(1-\kappa^{2})z}{1+i\kappa z}\\
 & =R(1-\kappa^{2})\sum_{n\geq1}(-i\kappa)^{n-1}z^{n}.
\end{align*}
The disk automorphism in the first expression verifies the image domain.
Summing the squares of the coefficients gives 
\[
\Var(\mu_{U_{\kappa}})=\frac{R^{2}(1-\kappa^{2})^{2}}{2}\sum_{n\geq1}\kappa^{2n-2}=\frac{R^{2}(1-\kappa^{2})}{2}.
\]
This expression takes every value in $(0,R^{2}/2]$ and tends to $0$
as $\kappa\uparrow1$. 
\end{proof}

\begin{remark} \label{op:rectangle} Small exit-law variance does
not imply small Gross area. Fix $a>0$ and set 
\[
R_{b}=\left\{ z\in\C:|\Re z|<\frac{a}{2b},\quad|\Im z|<\frac{b}{2}\right\} ,\qquad b>0.
\]
This bounded rectangle is symmetric about the real axis and vertically
convex. Its area is finite, so Proposition~\ref{op:geometric-equality},
applied with $s=1/2$, gives $U^{G}_{\mu_{R_{b}}}=R_{b}$. Consequently
\[
A(U^{G}_{\mu_{R_{b}}})=a,\qquad\rho_{1/2}(R_{b})=1.
\]
On the other hand, the real exit location lies in $[-a/(2b),a/(2b)]$,
and therefore 
\[
0<\Var(\mu_{R_{b}})\leq\frac{a^{2}}{4b^{2}}\longrightarrow0.
\]
Equivalently, the Gross boundary data tend to zero in normalized $L^{2}(\T)$,
while their spectral $H^{1/2}$ energies remain fixed. Thus Theorem~\ref{op:variance-range}
alone does not establish any limiting value for the area ratio of
its translated disks. \end{remark}

\section{An obstruction above first order}

\label{ce:section}

The contraction in Theorem~\ref{sp:domain} does not extend to all
higher orders, even when the competing domain has analytic boundary
and the exit law has compact support.

\begin{proposition}\label{ce:counterexample} Let 
\[
f(z)=z+\frac{1}{3}z^{2},\qquad U=f(\D),\qquad\mu=\mathcal{L}\!\left(\cos\Theta+\frac{1}{3}\cos2\Theta\right),
\]
where $\Theta$ is uniform on $(-\pi,\pi)$. Then $U$ is a bounded
simply connected $\mu$-domain containing $0$, and $\mu$ is centered,
nondegenerate, and compactly supported. Its Gross domain satisfies
\begin{equation}
\mathcal{E}_{s}(U^{G}_{\mu})<\infty\quad(0\leq s<3/2),\qquad\mathcal{E}_{s}(U^{G}_{\mu})=\infty\quad(s\geq3/2),\label{ce:threshold}
\end{equation}
whereas 
\begin{equation}
\mathcal{E}_{s}(U)=1+\frac{2^{2s}}{9}<\infty\qquad(s\geq0).\label{ce:polyenergy}
\end{equation}
Moreover, there exists $s_{*}\in(1,3/2)$ such that $\mathcal{E}_{s}(U^{G}_{\mu})>\mathcal{E}_{s}(U)$
for every $s\geq s_{*}$. For $s_{*}\leq s<3/2$, both sides are finite.
\end{proposition} 
\begin{proof}
The factorization 
\[
f(z)-f(w)=(z-w)\left(1+\frac{z+w}{3}\right)
\]
proves injectivity on a neighborhood of $\overline{\D}$: one may
take any disk of radius strictly between $1$ and $3/2$. Consequently
$U$ is bounded, simply connected, and has analytic Jordan boundary.
Since $f(0)=0$, conformal invariance identifies its real exit law
with the displayed law $\mu$. The zero integrals of the cosine modes
show that $\mu$ is centered; it is plainly nonconstant and compactly
supported. Equation~\eqref{ce:polyenergy} follows from the two nonzero
Taylor coefficients of $f$.

We compute the quantile $q$ of $\mu$. Put 
\[
P(x)=\frac{2}{3}x^{2}+x-\frac{1}{3}=\frac{2}{3}\left(x+\frac{3}{4}\right)^{2}-\frac{17}{24}.
\]
Thus $\cos\Theta+\frac{1}{3}\cos2\Theta=P(\cos\Theta)$. The minimum
is $t_{\min}=-17/24$ and the maximum is $t_{\max}=4/3$. At $t_{0}=-2/3$,
the roots of $P(x)=t_{0}$ are $-1$ and $-1/2$, so $F_{\mu}(t_{0})=1/3$.

For $t_{\min}\leq t\leq t_{0}$, let 
\[
a(t)=\frac{-3-\sqrt{17+24t}}{4},\qquad b(t)=\frac{-3+\sqrt{17+24t}}{4}.
\]
Both roots belong to $[-1,1]$, and the arcsine distribution of $\cos\Theta$
gives 
\[
F_{\mu}(t)=\frac{\arccos a(t)-\arccos b(t)}{\pi}.
\]
Set $r=F_{\mu}(t)$, $d=\pi r/2$, and $m=(\arccos a(t)+\arccos b(t))/2$.
Since $a+b=-3/2$, 
\[
\cos m=-\frac{3}{4\cos d},\qquad b-a=2\sin m\sin d=\frac{\sqrt{17+24t}}{2}.
\]
Squaring the latter identity and substituting the former expresses
$t$ as a function of $r$. For $t_{0}\leq t\leq t_{\max}$, only
the upper root bounds the event $P(\cos\Theta)\leq t$, and $F_{\mu}(t)=1-\arccos b(t)/\pi$.
Inverting the two formulas gives the exact quantile 
\begin{equation}
q(r)=\begin{cases}
-\dfrac{17}{24}+\dfrac{2}{3}\sin^{2}(\pi r/2)-\dfrac{3}{8}\tan^{2}(\pi r/2), & 0\leq r\leq1/3,\\[4pt]
\dfrac{2}{3}\cos^{2}(\pi r)-\cos(\pi r)-\dfrac{1}{3}, & 1/3\leq r\leq1.
\end{cases}\label{ce:quantile}
\end{equation}
Endpoint values are understood by continuity. These two analytic expressions
agree at $r=1/3$, but their derivatives satisfy 
\begin{equation}
q'(1/3-)=0,\qquad q'(1/3+)=\frac{\pi}{2\sqrt{3}}.\label{ce:jump}
\end{equation}
Also $q'(0)=q'(1)=0$. Let $\phi(\theta)=q(|\theta|/\pi)$ be the
increasing Gross datum, extended periodically, and write 
\[
\phi(\theta)=\sum_{n\geq1}a_{n}\cos(n\theta).
\]
 The constant term vanishes because $\mu$ is centered. Formula~\eqref{ce:quantile}
shows that $\phi$ is continuous and piecewise analytic, and is smooth
through $0$ and the periodic endpoints. Its first derivative has
jumps of size $J=1/(2\sqrt{3})$ at both $\theta=\pi/3$ and $\theta=-\pi/3$,
by \eqref{ce:jump}. Its distributional second derivative therefore
has the form 
\[
\phi''=g+J\delta_{\pi/3}+J\delta_{-\pi/3},
\]
where $g$ is piecewise smooth and has bounded variation. Twice integrating
by parts in the distributional sense yields 
\begin{align}
a_{n} & =-\frac{1}{\pi n^{2}}\left(\int^{\pi}_{-\pi}g(\theta)\cos(n\theta)\,d\theta+2J\cos(n\pi/3)\right)\nonumber \\
 & =-\frac{\cos(n\pi/3)}{\sqrt{3}\,\pi n^{2}}+O(n^{-3}).\label{ce:coefficients}
\end{align}
The error estimate follows by integration by parts on the finitely
many smooth pieces of $g$. In particular, $|a_{n}|=O(n^{-2})$, proving
convergence of $\sum n^{2s}|a_{n}|^{2}$ for $s<3/2$. Since $|\cos(n\pi/3)|\geq1/2$
for every integer $n$, \eqref{ce:coefficients} also gives $|a_{n}|\geq(4\sqrt{3}\,\pi)^{-1}n^{-2}$
for all sufficiently large $n$. Thus $\sum n^{2s}|a_{n}|^{2}$ diverges
for $s\geq3/2$. This proves \eqref{ce:threshold}; equivalently,
$\phi\in H^{s}(\T)$ for $s<3/2$ but $\phi\notin H^{3/2}(\T)$.

Finally, monotone convergence in the nonnegative coefficient sum implies
\[
\lim_{s\uparrow3/2}\mathcal{E}_{s}(U^{G}_{\mu})=\mathcal{E}_{3/2}(U^{G}_{\mu})=\infty.
\]
On the other hand, $\mathcal{E}_{s}(U)\leq17/9$ for $s\leq3/2$.
Choose $s_{*}\in(1,3/2)$ such that $\mathcal{E}_{s_{*}}(U^{G}_{\mu})>17/9$.
Monotonicity in $s$ then proves the strict reverse inequality for
$s_{*}\leq s<3/2$, with finite energies by \eqref{ce:threshold}.
For $s\geq3/2$, it follows directly from \eqref{ce:threshold} and
\eqref{ce:polyenergy}. 
\end{proof}

\section{Uniform upper energies and remaining questions}

\label{op:questions-section}

Let $\mu$ now be a centered, nondegenerate probability measure with
finite second moment. To discuss uniform upper bounds for the energies,
use the class 
\[
\mathcal{C}_{2}(\mu)=\left\{ U\subset\C:\begin{array}{l}
U\text{ is simply connected},\ 0\in U,\ \mu_{U}=\mu,\\
f_{U}\in H^{2}
\end{array}\right\} .
\]
This class is nonempty by the Gross construction. All its elements
satisfy $\mathcal{E}_{0}(U)=2\Var(\mu)$ by \eqref{op:h2-identities}.
Define 
\begin{equation}
\sigma_{\mu}=\sup\left(\{0\}\cup\left\{ s>0:\sup_{U\in\mathcal{C}_{2}(\mu)}\mathcal{E}_{s}(U)<\infty\right\} \right).\label{op:index}
\end{equation}
Including $0$ specifies the value when no positive exponent has a
uniform bound.

\begin{theorem} \label{op:index-range} The index in \eqref{op:index}
satisfies $0\leq\sigma_{\mu}\leq1/2$. \end{theorem} 
\begin{proof}
The lower bound follows from the definition. The following argument
proves directly the infinite-area consequence of the alternative construction
in~\cite{BMNew}. We construct a domain $V\in\mathcal{C}_{2}(\mu)$
with infinite area. Put $u(\theta)=q_{\mu}(\theta/(2\pi))$ for $0<\theta<2\pi$,
and let 
\[
F(z)=\int^{2\pi}_{0}\frac{e^{i\theta}+z}{e^{i\theta}-z}u(\theta)\,\frac{d\theta}{2\pi}.
\]
Then $F(0)=0$, $\Re F^{*}=u$, and $\|F\|^{2}_{\mathbf{H}^{2}}=2\|u\|^{2}_{L^{2}(dm)}>0$
by Parseval.

First suppose that $u$ is bounded and nonconstant. Its Stieltjes
measure $du$ on $(0,2\pi)$ is positive and nonzero. Integration
by parts, with the one-sided endpoint values, gives 
\[
-i(1-z)^{2}F'(z)=\frac{1}{\pi}\int_{(0,2\pi)}\frac{(1-z)(e^{i\theta}-1)}{e^{i\theta}-z}\,du(\theta).
\]
The identity 
\[
\Re\frac{(1-z)(e^{i\theta}-1)}{e^{i\theta}-z}=\frac{(1-\cos\theta)(1-|z|^{2})}{|e^{i\theta}-z|^{2}}
\]
shows that the real part of $-i(1-z)^{2}F'(z)$ is strictly positive.
For $\Re w>0$, set $\Psi(w)=-iF((w-1)/(w+1))$. With $z=(w-1)/(w+1)$,
we have $\Psi'(w)=-i(1-z)^{2}F'(z)/2$, so $\Re\Psi'(w)>0$. For distinct
$w_{1},w_{2}$ in the right half-plane, 
\[
\frac{\Psi(w_{2})-\Psi(w_{1})}{w_{2}-w_{1}}=\int^{1}_{0}\Psi'(w_{1}+t(w_{2}-w_{1}))\,dt
\]
has positive real part. Thus $\Psi$, and hence $F$, is univalent.

For general $u\in L^{2}$, let 
\[
u_{m}=\max(-m,\min(u,m))-\int\max(-m,\min(u,m))\,dm.
\]
These functions are bounded, nondecreasing and, for all sufficiently
large $m$, nonconstant. Their Schwarz integrals $F_{m}$ are univalent
by the preceding argument and converge to $F$ in $\mathbf{H}^{2}$,
hence locally uniformly. The nonconstant limit theorem for univalent
maps makes $F$ univalent. Its image $V=F(\D)$ is proper: otherwise
$F^{-1}$ would be a bounded nonconstant entire function, contrary
to Liouville's theorem. Lemma~\ref{pr:exit-time} gives $\E\tau_{V}<\infty$.
Lemma~\ref{pr:boundary-law} now gives $\mu_{V}=\mu$, so $V\in\mathcal{C}_{2}(\mu)$.

It remains to prove that $V$ has infinite area. Since the quantile
is nonconstant, there are $0<a<b<1$ and $\delta>0$ such that $q_{\mu}(b)-q_{\mu}(a)\ge\delta$.
For sufficiently small $\varepsilon>0$, 
\[
|u(x)-u(2\pi-y)|\ge\delta,\qquad0<x,y<\varepsilon.
\]
As $|e^{ix}-e^{-iy}|\le x+y$, it follows that 
\[
\iint_{\T^{2}}\frac{|u(\theta)-u(\varphi)|^{2}}{|e^{i\theta}-e^{i\varphi}|^{2}}\,dm(\theta)\,dm(\varphi)\ge\frac{\delta^{2}}{4\pi^{2}}\int^{\varepsilon}_{0}\int^{\varepsilon}_{0}\frac{dx\,dy}{(x+y)^{2}}=\infty.
\]
Writing $F(z)=\sum_{n\ge1}c_{n}z^{n}$, Parseval followed by Tonelli
gives the extended identity 
\[
\iint_{\T^{2}}\frac{|u(\theta)-u(\varphi)|^{2}}{|e^{i\theta}-e^{i\varphi}|^{2}}\,dm(\theta)\,dm(\varphi)=\sum_{k\in\mathbb{Z}}|k|\,|\widehat{u}(k)|^{2}=\frac{1}{2}\sum_{n\ge1}n|c_{n}|^{2}.
\]
For the first equality one uses $\int_{\T}|1-e^{ik\theta}|^{2}/|1-e^{i\theta}|^{2}\,dm(\theta)=|k|$;
this follows by expanding the finite geometric sum. Consequently $\mathcal{E}_{1/2}(V)=\infty$,
equivalently $\operatorname{Area}(V)=\infty$. Monotonicity of the
coefficient weights gives $\mathcal{E}_{s}(V)=\infty$ for every $s\ge1/2$.
No such exponent is uniformly bounded on $\mathcal{C}_{2}(\mu)$,
which proves $\sigma_{\mu}\le1/2$. 
\end{proof}

\begin{theorem}[Characterization of the upper-energy index] \label{sur:index-characterization}
Let $\mu\in\mathcal{M}_{2}$, let $X$ have law $\mu$, and let $h_{\mu}(t)=Q_{|X|}(1-t)$,
$0<t<1$, be the decreasing quantile of $|X|$, defined up to null
sets. For $0<s<1/2$, the following are equivalent: 
\begin{equation}
\sup_{U\in\mathcal{C}_{2}(\mu)}\mathcal{E}_{s}(U)<\infty\quad\Longleftrightarrow\quad\mathcal{E}_{s}(U^{G}_{\mu})<\infty\quad\Longleftrightarrow\quad L_{s}(\mu):=\int^{1}_{0}h_{\mu}(t)^{2}t^{-2s}\,dt<\infty.\label{sur:index-equivalences}
\end{equation}
Moreover, put 
\[
p_{*}(\mu)=\sup\{p\ge2:\E|X|^{p}<\infty\}\in[2,\infty].
\]
Then, with $1/\infty=0$, 
\begin{equation}
\boxed{{\displaystyle \sigma_{\mu}=\frac{1}{2}-\frac{1}{p_{*}(\mu)}.}}\label{sur:index-moments}
\end{equation}
In particular, $\sigma_{\mu}>0$ if and only if $\mu$ has a finite
moment of some order strictly greater than $2$. \end{theorem} 
\begin{proof}
The Gross domain belongs to $\mathcal{C}_{2}(\mu)$, so uniform finiteness
implies finiteness of its energy. If the latter is finite, Lemma~\ref{sur:hardy},
applied to its real boundary datum, gives $L_{s}(\mu)<\infty$. Conversely,
every $U\in\mathcal{C}_{2}(\mu)$ has an $\mathbf{H}^{2}$ conformal
map whose real boundary values have law $\mu$. Lemma~\ref{sur:univalent-reverse}
therefore gives the uniform estimate 
\[
\mathcal{E}_{s}(U)\le A_{s}L_{s}(\mu),\qquad U\in\mathcal{C}_{2}(\mu),
\]
proving \eqref{sur:index-equivalences}.

To identify the supremum of these exponents, set $p_{s}=2/(1-2s)$.
If $L_{s}(\mu)<\infty$, monotonicity of $h_{\mu}$ gives 
\[
L_{s}(\mu)\ge h_{\mu}(t)^{2}\int^{t}_{0}v^{-2s}\,dv=\frac{h_{\mu}(t)^{2}t^{1-2s}}{1-2s}.
\]
Hence $h_{\mu}(t)\le\sqrt{(1-2s)L_{s}(\mu)}\,t^{-1/p_{s}}$, and consequently
$\E|X|^{q}=\int^{1}_{0}h_{\mu}(t)^{q}\,dt<\infty$ for every $0<q<p_{s}$.
In the other direction, if $\E|X|^{q}<\infty$ for some $q>p_{s}$,
then $h_{\mu}(t)\le\|X\|_{L^{q}}t^{-1/q}$, so 
\[
L_{s}(\mu)\le\|X\|^{2}_{L^{q}}\int^{1}_{0}t^{-2s-2/q}\,dt<\infty.
\]
Thus every $s<1/2-1/p_{*}(\mu)$ with $s>0$ satisfies the equivalent
conditions in \eqref{sur:index-equivalences}, and no $s>1/2-1/p_{*}(\mu)$
with $s<1/2$ does. Proposition~\ref{op:index-range} excludes every
exponent $s\ge1/2$ from the set of uniformly bounded exponents. Taking
its supremum, with $0$ included as in \eqref{op:index}, proves \eqref{sur:index-moments},
including the cases $p_{*}(\mu)=2$ and $p_{*}(\mu)=\infty$. 
\end{proof}

The critical moment exponent $p_{*}(\mu)$ admits a classical interpretation
in terms of Hardy spaces and Brownian exit times. For every $U\in\mathcal{C}_{2}(\mu)$,
boundedness of harmonic conjugation on $L^{p}$ and Burkholder’s theorem
give 
\[
p_{*}(\mu)=h(U)=2bh(U),
\]
where $h(U)$ denotes the Hardy number and $bh(U)=\sup\{q>0:\mathbb{E}[\tau^{q}_{U}]<\infty\}$
is the Brownian--Hardy number. See \cite{GMCdonald,MG comb dom}.

\begin{example}{[}Heavy-tailed laws and the upper-energy index{]}

The upper-energy index can take any value in $(0,1/2)$, as illustrated
by the following families. All laws below are symmetric and have finite
variance.

$\smallskip$

$\noindent\textbf{Symmetric Pareto laws.}$

For $\alpha>2$, consider the density
\[
f_{\alpha}(x)=\frac{\alpha}{2}|x|^{-1-\alpha}\mathbf{1}_{\{|x|\ge1\}}.
\]
Since $\mathbb{E}|X|^{p}=\alpha/(\alpha-p)$ for $0\le p<\alpha$,
whereas the moment diverges for $p\ge\alpha$, we obtain
\[
p_{*}(\mu)=\alpha,\qquad\sigma_{\mu}=\frac{1}{2}-\frac{1}{\alpha}.
\]
More explicitly, $h_{\mu}(t)=t^{-1/\alpha}$, and hence
\[
L_{s}(\mu)=\begin{cases}
(1-2s-2/\alpha)^{-1}, & 0<s<\frac{1}{2}-\frac{1}{\alpha},\\[3pt]
\infty, & \frac{1}{2}-\frac{1}{\alpha}\le s<\frac{1}{2}.
\end{cases}
\]
Thus every index $\beta\in(0,1/2)$ is realized by choosing $\alpha=2/(1-2\beta)$.

$\smallskip$

$\noindent\textbf{Student's distributions.}$

For Student's \$t\$ distribution with \$\textbackslash nu>2\$ degrees
of freedom, the density is
\[
f_{\nu}(x)=\frac{\Gamma((\nu+1)/2)}{\sqrt{\nu\pi}\,\Gamma(\nu/2)}\left(1+\frac{x^{2}}{\nu}\right)^{-(\nu+1)/2}.
\]
Its density is asymptotic to a positive constant times $|x|^{-1-\nu}$
as $|x|\to\infty$. Consequently,
\[
p_{*}(\mu)=\nu,\qquad\sigma_{\mu}=\frac{1}{2}-\frac{1}{\nu}.
\]

$\smallskip$

$\noindent\textbf{Symmetric beta-prime laws.}$

For $a>0$ and $b>2$, let
\[
f_{a,b}(x)=\frac{|x|^{a-1}}{2B(a,b)(1+|x|)^{a+b}},\qquad x\ne0,
\]
where $B$ denotes the beta function. Direct integration gives
\[
\mathbb{E}|X|^{p}=\frac{B(a+p,b-p)}{B(a,b)}\quad(0\le p<b),
\]
and the moment is infinite for $p\ge b$. Therefore,
\[
p_{*}(\mu)=b,\qquad\sigma_{\mu}=\frac{1}{2}-\frac{1}{b}.
\]
In particular, the parameter $a$ changes the behavior of the density
near the origin without changing the upper-energy index.

$\smallskip$

$\noindent\textbf{Symmetrized inverse-gamma laws.}$

For $a>2$ and $\lambda>0$, consider

\[
f_{a,\lambda}(x)=\frac{\lambda^{a}}{2\Gamma(a)}|x|^{-a-1}e^{-\lambda/|x|},\qquad x\ne0.
\]
The substitution $y=\lambda/x$ yields
\[
\mathbb{E}|X|^{p}=\lambda^{p}\frac{\Gamma(a-p)}{\Gamma(a)}\quad(0\le p<a),
\]
with divergence for $p\ge a$. Hence

\[
p_{*}(\mu)=a,\qquad\sigma_{\mu}=\frac{1}{2}-\frac{1}{a}.
\]

\end{example}

\begin{example}[Different behavior at the critical index]

The value of $\sigma_{\mu}$ alone does not determine whether the
energies are uniformly bounded at $s=\sigma_{\mu}$. To illustrate
this distinction, fix $\alpha>2$ and $\beta\in\mathbb{R}$, and consider
the probability density
\[
f_{\alpha,\beta}(x)=C_{\alpha,\beta}|x|^{-1-\alpha}(\log|x|)^{-\beta}\mathbf{1}_{\{|x|\ge e\}},
\]
where $C_{\alpha,\beta}>0$ is the normalizing constant. For every
$\beta$, the moments of order $p<\alpha$ are finite, whereas those
of order $p>\alpha$ are infinite. Thus
\[
p_{*}(\mu)=\alpha,\qquad s_{c}:=\sigma_{\mu}=\frac{1}{2}-\frac{1}{\alpha}.
\]
We now examine the critical exponent. Integration of the density gives
\[
\mathbb{P}(|X|>x)\sim\frac{2C_{\alpha,\beta}}{\alpha}x^{-\alpha}(\log x)^{-\beta}\qquad(x\to\infty).
\]
Inverting this tail estimate yields
\[
h_{\mu}(t)\asymp t^{-1/\alpha}\bigl(\log(1/t)\bigr)^{-\beta/\alpha}\qquad(t\downarrow0).
\]
Here $\asymp$ means that the ratio is bounded above and below by
positive constants. Since $2s_{c}=1-2/\alpha$, the integrand defining
$L_{s_{c}}(\mu)$ is comparable near zero to
\[
\frac{1}{t\bigl(\log(1/t)\bigr)^{2\beta/\alpha}}.
\]
It follows that
\[
L_{s_{c}}(\mu)<\infty\quad\Longleftrightarrow\quad\beta>\frac{\alpha}{2}.
\]
By the characterization theorem, therefore,
\[
\sup_{U\in\mathcal{C}_{2}(\mu)}\mathcal{E}_{s_{c}}(U)<\infty\quad\Longleftrightarrow\quad\beta>\frac{\alpha}{2}.
\]
Thus laws with the same moment threshold and the same upper-energy
index can exhibit different behavior at the critical exponent. The
Lorentz criterion detects this

distinction, which is not captured by the index formula alone.

\end{example}

\begin{remark} Formula~\eqref{sur:index-moments} determines the
index, not membership at its critical exponent. If $0<\sigma_{\mu}<1/2$,
uniform finiteness at $s=\sigma_{\mu}$ is decided by the integral
in \eqref{sur:index-equivalences}; no endpoint moment assumption
is implicit in the formula. When $\sigma_{\mu}=1/2$, the uniform
bound fails at that endpoint by Proposition~\ref{op:index-range}.
\end{remark}

\begin{corollary}[Attainment of the maximal energy below the index]
\label{sur:maximum-attainment} Let $\mu\in\mathcal{M}_{2}$ and $0<s<\sigma_{\mu}$.
Then there exists $U_{\max}\in\mathcal{C}_{2}(\mu)$ such that 
\[
\mathcal{E}_{s}(U_{\max})=\max_{U\in\mathcal{C}_{2}(\mu)}\mathcal{E}_{s}(U)<\infty.
\]
\end{corollary} 
\begin{proof}
Choose $s<r<\sigma_{\mu}$. By Corollary~\ref{sur:index-characterization},
there is a finite constant $M$ such that $\mathcal{E}_{r}(U)\le M$
for every $U\in\mathcal{C}_{2}(\mu)$. Consequently the normalized
conformal maps of this entire class form exactly the set $\mathcal{D}_{2r,M}(\mu)$
from Proposition~\ref{tr:compactness}: conversely, finite $\mathcal{E}_{r}$
implies $\mathbf{H}^{2}$ integrability. That proposition makes this
nonempty set compact in $\mathcal{D}_{2s}$. Since $\mathcal{E}_{s}(U)=\|f_{U}\|^{2}_{\mathcal{D}_{2s}}$
is continuous in this norm, it attains its maximum there. 
\end{proof}

For finite-energy competitors, the ratio question is settled at the
area exponent and throughout $1/2\le s\le1$ by Theorem~\ref{sur:log-surjectivity}.
Theorem~\ref{sur:positive-infimum} rules out surjectivity onto $(0,1]$
below $1/2$. The following questions remain.
\begin{enumerate}
\item Does $\mathcal{P}(\mu)\subsetneq\mathcal{B}(\mu)$ hold for every
$\mu\in\mathcal{M}_{*}$? Theorem~\ref{pr:strict-inclusion} establishes
strict inclusion for a particular compactly supported law. 
\item For $0<s<1/2$, determine the sharp constant 
\[
\kappa^{\mathrm{opt}}_{s}=\inf\{\rho_{s}(U):0<\mathcal{E}_{s}(U)<\infty\}.
\]
Theorem~\ref{sur:positive-infimum} proves that it is positive. Is
the infimum attained, and what are the extremal domains? 
\item If $0<\sigma_{\mu}<1/2$ and $L_{\sigma_{\mu}}(\mu)<\infty$, is the
maximal energy attained at the critical exponent $\sigma_{\mu}$?
Corollary~\ref{sur:maximum-attainment} proves attainment below that
exponent; its compactness argument uses a strict gap between exponents. 
\item Determine which exponents above $1$, or which additional hypotheses
on the quantile, permit the minimization in Theorem~\ref{sp:domain}.
Proposition~\ref{ce:counterexample} rules out a universal extension
for exponents sufficiently close to $3/2$ and for every exponent
at least $3/2$. 
\end{enumerate}

\section{A computational perspective}

The exit-law formulation suggests estimating the real exit law from
Brownian simulations and reconstructing a regularized Gross map from
the estimated quantile. Developing and analyzing such a scheme remains
an open computational problem, requiring control of exit detection,
sampling error, and reconstruction error. The area criterion above
identifies an obstruction to using raw empirical quantiles: every
nondegenerate finitely supported law has an internal support gap and
therefore infinite Gross area. Consequently, convergence of empirical
laws in $W_{2}$ alone does not guarantee convergence of minimal areas.
A finite-area approximation therefore requires regularization whose
effect on the target law and its logarithmic energy is controlled.

\paragraph{Data availability.}

No empirical datasets were generated or analyzed. All examples and
proofs are analytical.

\raggedright


\begin{thebibliography}{10}
\bibitem{Gross} R. Gross, A conformal Skorokhod embedding, \emph{Electronic
Communications in Probability} \textbf{24} (2019), paper 68, 11 pp.
\href{https://doi.org/10.1214/19-ECP272}{doi:10.1214/19-ECP272}.

\bibitem{BMRemarks} M. Boudabra and G. Markowsky, Remarks on Gross'
technique for obtaining a conformal Skorohod embedding of planar Brownian
motion, \emph{Electronic Communications in Probability} \textbf{25}
(2020), paper 20, 13 pp. \href{https://doi.org/10.1214/20-ECP300}{doi:10.1214/20-ECP300}.

\bibitem{BMNew} M. Boudabra and G. Markowsky, A new solution to the
conformal Skorokhod embedding problem and applications to the Dirichlet
eigenvalue problem, \emph{Journal of Mathematical Analysis and Applications}
\textbf{491}(2) (2020), article 124351. \href{https://doi.org/10.1016/j.jmaa.2020.124351}{doi:10.1016/j.jmaa.2020.124351}.

\bibitem{BoudabraEndpoint} M. Boudabra, A note on the planar Skorokhod
embedding problem, \emph{Bulletin of the Australian Mathematical Society}
(2026), published online 12 February 2026, 7 pp. \href{https://doi.org/10.1017/S0004972726100987}{doi:10.1017/S0004972726100987}.

\bibitem{Energy} M. Becher, M. Boudabra, and F. Haggui, Skorokhod
energy of planar domains, \emph{Journal of Mathematical Analysis and
Applications} \textbf{551}(1) (2025), article 129641. \href{https://doi.org/10.1016/j.jmaa.2025.129641}{doi:10.1016/j.jmaa.2025.129641}.

\bibitem{Baernstein} A. Baernstein II, Some sharp inequalities for
conjugate functions, \emph{Indiana University Mathematics Journal}
\textbf{27}(5) (1978), 833--852. \href{https://doi.org/10.1512/iumj.1978.27.27055}{doi:10.1512/iumj.1978.27.27055}.

\bibitem{Pruss} A. R. Pruss, Steiner symmetry, horizontal convexity,
and some Dirichlet problems, \emph{Indiana University Mathematics
Journal} \textbf{46}(3) (1997), 863--895.

\bibitem{RoncalStinga} L. Roncal and P. R. Stinga, Fractional Laplacian
on the torus, \emph{Communications in Contemporary Mathematics} \textbf{18}
(2016), article 1550033, 26 pp. \href{https://doi.org/10.1142/S0219199715500339}{doi:10.1142/S0219199715500339}.

\bibitem{CCM} X. Cabré, G. Csató, and A. Mas, Periodic solutions
to integro-differential equations: variational formulation, symmetry,
and regularity, \emph{Communications in Partial Differential Equations}
\textbf{50}(1--2) (2025), 162--210. \href{https://doi.org/10.1080/03605302.2024.2441851}{doi:10.1080/03605302.2024.2441851}.

\bibitem{StrictRearrangement} G. Csató and A. Mas, Strict rearrangement
inequalities: Nonexpansivity and periodic Gagliardo seminorms, \emph{Transactions
of the American Mathematical Society} \textbf{378}(10) (2025), 7163--7197.
\href{https://doi.org/10.1090/tran/9510}{doi:10.1090/tran/9510}.

\bibitem{BrownianConvergence} Boudabra, M., Becher, M. \& Haggui,
F. Brownian convergence of planar domains and stability of the planar
Skorokhod embedding problem. Arch. Math. 127, 95--111 (2026). 

\bibitem{BoundedGross} M. Boudabra, D. Kroumi, and B. Smii, On the
boundedness of Gross' solution to the planar Skorokhod embedding problem,
arXiv:2412.02913 (2024). \href{https://arxiv.org/abs/2412.02913}{arXiv:2412.02913}.

\bibitem{Burkholder} D. L. Burkholder, Exit times of Brownian motion,
harmonic majorization, and Hardy spaces, \emph{Advances in Mathematics}
\textbf{26}(2) (1977), 182--205. \href{https://doi.org/10.1016/0001-8708(77)90029-9}{doi:10.1016/0001-8708(77)90029-9}.

\bibitem{DurenHardy} P. L. Duren, \emph{Theory of $\mathbf{H}^{p}$
Spaces}, Pure and Applied Mathematics, vol.~38, Academic Press, New
York, 1970.

\bibitem{MortersPeres} P. Mörters and Y. Peres, \emph{Brownian Motion},
Cambridge Series in Statistical and Probabilistic Mathematics, vol.~30,
Cambridge University Press, Cambridge, 2010. \href{https://doi.org/10.1017/CBO9780511750489}{doi:10.1017/CBO9780511750489}.

\bibitem{DurenUnivalent} P. L. Duren, \emph{Univalent Functions},
Grundlehren der mathematischen Wissenschaften, vol.~259, Springer-Verlag,
New York, 1983.

\bibitem{Santambrogio} F. Santambrogio, \emph{Optimal Transport for
Applied Mathematicians: Calculus of Variations, PDEs, and Modeling},
Progress in Nonlinear Differential Equations and Their Applications,
vol.~87, Birkhäuser, Cham, 2015. \href{https://doi.org/10.1007/978-3-319-20828-2}{doi:10.1007/978-3-319-20828-2}.

\bibitem{Peetre} J. Peetre, Espaces d'interpolation et théorème de
Soboleff, \emph{Annales de l'Institut Fourier} \textbf{16}(1) (1966),
279--317. \href{https://doi.org/10.5802/aif.232}{doi:10.5802/aif.232}.

\bibitem{QianShi} R. Qian and Y. Shi, Univalent functions in Dirichlet-type
spaces, \emph{Filomat} \textbf{30}(5) (2016), 1213--1218. \href{https://doi.org/10.2298/FIL1605213Q}{doi:10.2298/FIL1605213Q}.

\bibitem{BennettSharpley} C. Bennett and R. Sharpley, \emph{Interpolation
of Operators}, Academic Press, 1988.

\bibitem{Popoviciu} T. Popoviciu, Sur les %
\mbox{%
équations%
} algébriques ayant toutes leurs racines réelles, \emph{Mathematica
(Cluj)} \textbf{9} (1935), 129--145. \href{https://ictp.acad.ro/fr/sur-les-equations-algebriques-ayant-toutes-leurs-racines-reelles/}{Original paper and transcription}.

\bibitem{Rudin} W. Rudin, Real and Complex Analysis 3rd ed., McGraw-{}-Hill,
New York, 1987.

\bibitem{Carleson} L. Carleson, On convergence and growth of partial
sums of Fourier series, Acta Mathematica $\textbf{116}$(1966), 135-{}-157.

\bibitem{Kesavan} Kesavan, S. 2006. Symmetrization and Applications.
Vol. 3, Series in Analysis. Singapore: World Scientific Publishing
Company.

\bibitem{MG comb dom} M. Boudabra and G. Markowsky, “On the finiteness
of moments of the exit time of planar Brownian motion from comb domains,”
Annales Fennici Mathematici 46 (2021), no. 1, 527--536.

\bibitem{GMCdonald} G. Markowsky and C. McDonald, “A collection of
results relating the geometry of plane domains and the exit time of
planar Brownian motion, II,” Computational Methods and Function Theory
(2026), published online 12 May 2026.

\bibitem{Ohrem} S. Ohrem and W. Reichel, Existence of traveling breather
solutions to cubic nonlinear Maxwell equations in waveguide geometries,
Calculus of Variations and Partial Differential Equations (2026).

\bibitem{Ba=0000F1uelos} R. Bañuelos and T. Carroll, Brownian motion
and the fundamental frequency of a drum, Duke Mathematical Journal
75 (1994), no. 3, 575--602.

\end{thebibliography}
\end{document}